\documentclass[a4paper,10pt]{siamart190516}
\usepackage{amsfonts}
\usepackage{MnSymbol}
\usepackage{enumerate}
\usepackage{bbm}
\usepackage{xcolor}
\usepackage{algorithm}
\usepackage[noend]{algpseudocode}
\usepackage{subfigure}
\usepackage{hyperref}

\newsiamremark{remark}{Remark}
\newsiamremark{hypothesis}{Hypothesis}
\crefname{hypothesis}{Hypothesis}{Hypotheses}
\newsiamthm{claim}{Claim}

\headers{The AZ Algorithm}{V.~Copp\'e, D.~Huybrechs, R.~Matthysen, M.~Webb}

\newcommand{\vertiii}[1]{{\left\vert\kern-0.25ex\left\vert\kern-0.25ex\left\vert #1 
		\right\vert\kern-0.25ex\right\vert\kern-0.25ex\right\vert}}
	
\newcommand{\R}{\mathbb{R}}
\ifdefined\C
\renewcommand{\C}{\mathbb{C}}
\else
\newcommand{\C}{\mathbb{C}}
\fi

\newcommand{\Z}{\mathbb{Z}}

\newcommand{\rank}{\mathrm{rank}}
\newcommand{\diag}{\mathrm{diag}}

\newcommand{\eps}{\varepsilon}
\newcommand{\FF}{\mathrm{F}}
\newcommand{\ee}{\mathrm{e}}
\newcommand{\ii}{\mathrm{i}}

\title{The AZ algorithm for least squares systems with \\ a known incomplete generalized inverse}
\author{Vincent Copp\'e\thanks{Department of Computer Science, KU Leuven, Belgium. \texttt{vincent.coppe@cs.kuleuven.be}, \texttt{daan.huybrechs@cs.kuleuven.be}, \texttt{roelmatthysen@gmail.com}} \and Daan Huybrechs\footnotemark[1] \and Roel Matthysen\footnotemark[1] \and Marcus Webb\thanks{Department of Mathematics, University of Manchester, UK. \texttt{marcus.webb@manchester.ac.uk}. This author is grateful to FWO Research Foundation Flanders for the postdoctoral fellowship he enjoyed during the research for this paper.}}

\begin{document}

\maketitle

\begin{abstract}
We introduce an algorithm for the least squares solution of a rectangular linear system $Ax=b$, in which $A$ may be arbitrarily ill-conditioned. We assume that a complementary matrix $Z$ is known such that $A - AZ^*A$ is numerically low rank. Loosely speaking, $Z^*$ acts like a generalized inverse of $A$ up to a numerically low rank error. We give several examples of $(A,Z)$ combinations in function approximation, where we can achieve high-order approximations in a number of non-standard settings: the approximation of functions on domains with irregular shapes, weighted least squares problems with highly skewed weights, and the spectral approximation of functions with localized singularities. The algorithm is most efficient when $A$ and $Z^*$ have fast matrix-vector multiplication and when the numerical rank of $A - AZ^*A$ is small.
\end{abstract}


\begin{AMS}
  [MSC2010] Primary: 65F20, Secondary: 65F05, 68W20
\end{AMS}

% !TEX root = AZalgorithmV3.tex
\section{Introduction}
\label{s:intro}

The topic of this paper arose from the study of function approximation using frames \cite{adcock2019frames,adcock2018fna2}. Frames are sets of functions that, unlike a basis, may be redundant. This leads to improved flexibility in a variety of settings which we illustrate below, but it also comes at a cost: the approximation problem is equivalent to solving a linear system $Ax=b$ that is, owing to the redundancy of the frame, ill-conditioned. The central result of \cite{adcock2019frames,adcock2018fna2} is that the function at hand can nevertheless be approximated to high accuracy via a regularized least squares approximation. Thus, the linear system is rectangular, $A \in \mathbb{C}^{M \times N}$ with $M >  N$, and the cost of a direct solver scales as ${\mathcal O}(MN^2)$, where typically $M$ is at least linear in $N$.

The particular type of frame that has received the most attention in this context corresponds to Fourier extension or Fourier continuation \cite{boyd2002fourierextension,bruno2007continuation,huybrechs2010fourierextension}, where a smooth but non-periodic function on a domain $\Omega$ is approximated by a Fourier series on a larger bounding box. Here, one can think of the redundancy as corresponding to different extensions of $f$ from $\Omega$ to the bounding box. Fast algorithms for the construction of Fourier extension approximations were proposed by Lyon in \cite{lyon2011fastcontinuation} for univariate problems and by Matthysen and Huybrechs in \cite{matthysen2016fastfe,matthysen2017fastfe2d} for more general problems in one and two dimensions. These algorithms have the form of the proposed AZ-algorithm in the special case where $Z = A$.
%\footnote{It should be noted that the univariate Fourier extension problem is equivalent to the problem of bandlimited signal extrapolation, for which a popular algorithm is the Papoulis-Gerchberg algorithm \cite{gerchberg1974superresolution,papoulis1975bandlimited}. The latter corresponds precisely to conjugate gradients (CG) applied to the normal equations of the least squares problem \cite{strohmer1995extrapolation}. While this algorithm provides sufficient accuracy for the purposes of signal and image processing, it falls short of providing high accuracy to numerical approximation problems due to the ill-conditioning of $A$ and corresponding slow convergence of CG.}

Thus, the AZ algorithm is a generalization of the algorithms proposed for Fourier extension in \cite{lyon2011fastcontinuation,matthysen2016fastfe,matthysen2017fastfe2d}. In this paper we analyze the algorithm from the point of view of linear algebra, rather than approximation theory. We supply error estimates when using techniques from randomized linear algebra for the first step of the AZ algorithm \cite{halko2011finding}. In contrast, the previous analysis in \cite{matthysen2016fastfe,matthysen2017fastfe2d} has focused exclusively on studying the numerical rank of this system, in the specific setting of Fourier extension problems. We characterize the scope of the algorithm and provide several examples in different settings. In our analysis we focus on backward error, as the size of the residual is directly related to the approximation error in the examples given.

%{\color{red} TODO: discussion of standard numerical methods and citing Golub/Van Loan, Lawson and Hanson book on least squares, iterative methods, }
There is a rich body of literature in numerical linear algebra on solution methods for least squares problems. We refer the reader to the standard references \cite{Golub2013,lawson1996leastsquares}. We note that direct solvers typically exhibit cubic complexity in the dimension of the problem: the aim of the AZ algorithm is to reduce that complexity by exploiting the specific structure of certain least squares problems. This leads to an algorithm with several steps, in which existing methods for least squares problems can be used for systems with lower rank. In particular we study SVD and QR based methods for step 1 of the algorithm. We briefly comment on the applicability of iterative solvers for least squares problems, such as LSQR \cite{Paige1982} and LSMR \cite{Fong2010}, near the end of the paper.

The structure of the paper is as follows. We formulate the algorithm in \S\ref{s:AZ} and state some general algebraic properties. We analyze the backward error of solvers for low-rank systems based on randomized SVD or randomized QR in \S\ref{s:lowrank}. After the analysis, we illustrate the algorithms with a sequence of examples: Fourier extension approximation in \S\ref{s:fourext}, approximation using weighted bases in \S\ref{s:sumframes}, spline-based extension problems in \S\ref{s:splines} and weighted least squares problems in \S\ref{s:weighted}. We end the paper with some concluding remarks in \S\ref{s:conclusion}.

% !TEX root = AZalgorithmV3.tex
\section{The AZ algorithm} \label{s:AZ}

Given a linear system $Ax=b$, with $A \in \mathbb{C}^{M \times N}$, and an additional matrix $Z \in \mathbb{C}^{M \times N}$, the AZ algorithm consists of three steps.

\begin{algorithm}
  \caption{The AZ algorithm}\label{alg:AZ}
    {\bf Input:} $A,Z \in \C^{M\times N}$, $b\in\C^M$ \\
     {\bf Output:} $x\in\C^N$ such that $Ax \approx b$ in a least squares sense
  \begin{algorithmic}[1]
    \State Solve $(I-AZ^*)Ax_1 = (I-AZ^*)b$
    \State $x_2 \gets Z^*(b-Ax_1)$
    \State $x \gets x_1 + x_2$
  \end{algorithmic}
\end{algorithm}

The algorithm does not specify which matrix $Z$ to use, nor which solver to use in step 1. The intention is to choose $Z$ such that the matrix $(I-AZ^*)A$ approximately has low rank. Since $(I-AZ^*)A = A - AZ^*A$, this property can be thought of as $Z^*$ being an incomplete generalized inverse of $A$. We analyze a few choices of solvers in \S\ref{s:lowrank} and we remark on the choice of $Z$ in later sections. Remarkably, the following two statements are independent of these choices.

\begin{lemma}[AZ Lemma]\label{lem:AZ}
  Let $\hat{x} = \hat{x}_1+\hat{x}_2$ be output from the AZ algorithm. Then the final residual is equal to the residual of step 1.
  \begin{proof}
    A simple expansion of the final residual yields
    \begin{eqnarray*}
    b - A\hat{x} &=& b - A\hat{x}_1 - A\hat{x}_2\\
    &=& b - A\hat{x}_1 - AZ^*(b-A\hat{x}_1) \\
    &=& (I-AZ^*)(b-A\hat{x}_1).
    \end{eqnarray*}
This is precisely the residual from step 1 of the algorithm.
    \end{proof}
  \end{lemma}

Thus, the accuracy with which step 1 is solved determines the accuracy of the overall algorithm, at least when accuracy is measured in terms of the size of the residual. Furthermore, the computational cost of the algorithm is also shifted to step 1, because step 2 and 3 only involve matrix-vector multiplication and addition.

In the following statement, we use the terminology of a \emph{stable least squares fit}. A stable least squares fit corresponds to a solution vector that has itself a moderate norm and that yields a small residual. The relevance of a stable least squares fit in the setting of our examples later in this paper is that, if a stable least squares fit exists, then numerical methods exist that can reliably find such a solution, no matter how ill-conditioned the linear system is. 

One way to do compute a stable least squares fit is based on a regularized SVD and this case was studied for function approximation using frames in \cite{adcock2019frames}. The following lemma proves that the AZ algorithm can also, in principle, find a stable least squares fit (if it exists).

\begin{lemma}[Stable least squares fitting]
  Let $A \in \C^{M\times N}, b\in\C^M$, and suppose there exists $\tilde{x} \in \C^N$ such that
  \begin{equation}
  \|b-A\tilde{x}\|_2 \leq \tau, \qquad \|\tilde{x}\|_2 \leq C,
  \end{equation}
  for $\tau,C > 0$. Then there exists a solution $\hat{x}_1$ to step 1 of the AZ algorithm such that the residual of the computed vector $\hat{x} = \hat{x}_1 + \hat{x}_2$ satisfies,
  \begin{equation}
   \|b-A \hat{x}\|_2 \leq \|I-AZ^*\|_2 \tau, \qquad \|\hat{x}\|_2 \leq C + \|Z^*\|_2 \tau
  \end{equation}
  \begin{proof}
    Take $\hat{x}_1 = \tilde{x}$. Then Lemma \ref{lem:AZ} implies $b-A\hat{x} = (I-AZ^*)(b-A\tilde{x})$, which gives us the first inequality. The second inequality follows from $\hat{x} = \tilde{x} + Z^*(b-A\tilde{x})$.
    \end{proof}
  \end{lemma}

\subsection{Choosing the matrix $Z$}

In principle, the AZ algorithm can be used with any matrix $Z$. However, we must choose $Z$ wisely if we want efficiency and accuracy. If the norm of $Z^*$ is large then Lemma 2.2 suggests that we may fail to obtain a stable least squares fit. If $A - AZ^*A$ is not low rank then the cost of solving the least squares system in step 1 may be expensive. Ideally, both $A$ and $Z^*$ also have fast matrix-vector multiplies.

Two extreme choices are $Z = 0$ and $Z = \left(A^\dagger\right)^*$. In the case $Z=0$, step 1 solves the original system and step 2 returns zero. We do not gain any efficiency over simply solving the original system. In the latter case, it is exactly the opposite: if $Z=\left(A^\dagger\right)^*$, then the first step yields $x_1=0$ and the original system is solved in the second step using the pseudoinverse of $A$. This may be efficient, but we would need to know the pseudoinverse of $A$!

Any other choice of $Z$ leads to a situation in between: part of the solution is found using a direct solver in step 1, part of it is found by multiplication with $Z^*$ in step 2. We desire a matrix $Z$ such that $A - AZ^*A$ has low numerical rank, so that step 1 can be solved efficiently, and $Z^*$ must be readily available with an efficient matrix-vector multiply if possible. Much like a preconditioner, we choose $Z$ satisfying these properties using some a priori information about the underlying problem. 

The following lemma gives a general relationship between $A$ and $Z$ which would guarantee $A - AZ^* A$ to be numerically low rank.

%Compare to preconditioning.

\begin{lemma}\label{lem:splitting}
  Suppose that $A,Z\in\C^{M\times N}$ satisfy
  \begin{equation*}
  A = W + L_1 + E_1, \qquad Z^* = W^\dagger + L_2 + E_2,
  \end{equation*}
  where $\rank(L_1),\rank(L_2) \leq R$ and $\|E_1\|_{\FF},\|E_2\|_{\FF} \leq \eps$. Here $W^\dagger$ is the Moore-Penrose pseudoinverse. Then
  \begin{equation*}
  A - AZ^*A = L + E,
  \end{equation*}
  where $\rank(L) \leq 3R$ and
  \begin{equation}\label{eqn:Ebound}
  \|E\|_{\FF} \leq \eps \left(1 + \|I-AZ^*\|_2 + \|A\|_2^2\right) + \eps^2 \|A\|_2.
  \end{equation}
  The result is exactly the same if the norms on $E_1$, $E_2$ and $E$ are changed to $\|\cdot\|_2$.
  \begin{proof}
    We simply expand $A-AZ^*A$. For clarity we expand terms gradually.
    \begin{eqnarray*}
    A - AZ^*A &=& W + L_1 + E_1 - AZ^*W - AZ^*(L_1+E_1) \\
              &=& W + L_1 + E_1 - AW^\dagger W - A(L_2+E_2)W - AZ^*(L_1+E_1) \\
              &=& W + L_1 + E_1 - WW^\dagger W - (L_1+E_1)W^\dagger W\\
              & &\qquad -~A(L_2+E_2)W - AZ^*(L_1+E_1)
    \end{eqnarray*}
    Since $W^\dagger$ is a generalized inverse of $W$, we have $W-WW^\dagger W = 0$, so
    \begin{equation*}
    A - AZ^*A =  L_1 + E_1 - (L_1+E_1)W^\dagger W - A(L_2+E_2)W - AZ^*(L_1+E_1).
    \end{equation*}
    Now, writing $A(L_2+E_2)W = AL_2 W + AE_2A - AE_2L_1 - AE_2E_1$ and splitting up the low rank and small norm parts, we obtain,
    \begin{eqnarray*}
    A-AZ^*A &=& \underbrace{(I-AZ^* - AE_2 )L_1 - L_1W^\dagger W - AL_2 W}_{L} \\
            & & \qquad + \qquad \underbrace{ (I-AZ^*)E_1 - W W^\dagger E_1 - AE_2 A + AE_2E_1}_{E}
    \end{eqnarray*}
    It is clear that $\rank(L) \leq 3R$ because of the three occurrences of $L_1$ and $L_2$. The Frobenius norm of $E$ is bounded above by
    \begin{eqnarray*}
    \|E\|_{\FF} &\leq& \|(I-AZ^*)E_1\|_{\FF} + \|W W^\dagger E_1\|_{\FF} + \|AE_2 A\|_{\FF} + \|AE_2E_1\|_{\FF}\\
                     &\leq&\|E_1\|_{\FF}\|I-AZ^*\|_2 + \|WW^\dagger\|_2 \|E_1\|_{\FF} + \|E_2\|_{\FF}\|A\|_2^2 + \|E_1\|_F \|E_2\|_2\|A\|_2
    \end{eqnarray*}
    Here we have used the inequality $\|BC\|_{\FF} \leq \|B\|_2\|C\|_{\FF}$ for general rectangular matrices $B$ and $C$\footnote{The proof of this is a one-liner: $\|BC\|_{\FF}^2 = \sum_{j} \|Bc_j\|_2^2 \leq \|B\|_2^2 \sum_{j} \|c_j\|_2^2 = \|B\|_2^2\|C\|_{\FF}^2 $}.
  
    Because $W^\dagger$ is the Moore-Penrose pseudoinverse of $W$, the matrix $W W^\dagger$ is an orthogonal projection onto the image of $W$, so $\|WW^\dagger\|_2 \leq 1$. Furthermore, $\|E_2\|_2 \leq \|E_2\|_{\FF} \leq \eps$, which enables us to arrive at the desired bound on $\|E\|_{\FF}$ in equation \eqref{eqn:Ebound}.
  
  In order to prove the same result with the Frobenius norms on $E_1$, $E_2$ and $E$ replaced by $\|\cdot\|_2$, it is readily checked that the exact same proof holds with all instances of $\FF$ replaced by $2$. In particular, we have the inequality $\|BC\|_2 \leq \|B\|_2\|C\|_2$.
  \end{proof}
\end{lemma}

In our examples, $A$ and $Z$ are determined by analytical means, which are appli\-cation-specific. We do not explicitly compute the $W$ and $L_1,L_2$ matrices of the above lemma. The effectiveness of the algorithm merely relies on the fact that these matrices exist. It may be possible to compute a suitable $Z$ numerically using a partial SVD of $A$, but that is unlikely to be efficient.

\section{Fast randomized algorithms for numerically low-rank systems}\label{s:lowrank}

In this section we discuss fast algorithms for the solution of the system $A x = b$, where $A \in \C^{M\times N}$, with $M\geq N$, has \emph{epsilon rank} $r$. We write $\rank_\eps (A) = r$ and in this paper it means that there exists $L,E \in \C^{M\times N}$ such that
\begin{equation}\label{eqn:lowrank}
A = L + E, \text{ where } \rank(L) = r \text{ and } \|E\|_\FF \leq \eps.
\end{equation}
It is important to note that we have used the Frobenius norm here, which implies $\sum_{k > r} \sigma_k^2 \leq \eps^2$, where $\sigma_{r+1},\ldots,\sigma_{N}$ are the $N-r$ smallest singular values of $A$. This is necessary for our proofs of the error bounds. Throughout this section there are points where the Frobenius norm appears to have been used unnecessarily where the 2-norm could have been used, but in all cases we do not believe that the final bounds will be improved significantly by changing to the 2-norm. This is due to the fact that we are not aware of an effective version of Proposition \ref{prop:Gaussianbounds1} which bounds 2-norms of the relevant random matrices purely in terms of 2-norms of other matrices.

We make use of Gaussian random matrices, $\Omega \sim \mathcal{N}(0,1;\R^{n\times k})$, which are $n \times k$ matrices whose elements are independent Gaussian random variables with mean $0$ and variance $1$.
\begin{proposition}[\cite{halko2011finding}]\label{prop:Gaussianbounds1}
  Let $\Omega \sim \mathcal{N}(0,1;\R^{n\times k})$, and let $S,T \in \R^{n\times k}$ be fixed matrices. Then for all $u \geq 0$,
  \begin{equation*}
  \mathbb{E}\left\{ \|S^* \Omega  T\|_\FF \right\} \leq \|S\|_\FF\|T\|_\FF, \qquad \mathbb{P}\left\{ \|S^* \Omega  T\|_\FF \geq (1+u) \cdot \|S\|_\FF\|T\|_\FF  \right\} \leq \ee^{-\frac{u^2}{2}}
  \end{equation*}
\end{proposition}

\begin{proposition}[\cite{halko2011finding}]\label{prop:Gaussianbounds2}
  Let $\Omega \sim \mathcal{N}(0,1;\R^{r \times (r+p)})$ with $p \geq 4$. Then for all $s\geq1$,
  \begin{equation*}
  \mathbb{E}\left\{ \|\Omega^\dagger \|_{\FF} \right\} = \sqrt{\frac{r}{p-1}}, \qquad \mathbb{P}\left\{ \|\Omega^\dagger \|_{\FF} \geq  s\cdot \sqrt{\frac{3r}{p+1}}\right\} \leq s^{-p}
  \end{equation*}
\end{proposition}

Note that we intend to apply the results of this section to step 1 of the AZ algorithm. Thus, matrix $A$ in this section actually corresponds to matrix $A-AZ^*A$ in Algorithm~\ref{alg:AZ}.

\subsection{Truncated SVD solvers}

\begin{algorithm}[h!]
  \caption{Truncated SVD solver \cite{Golub2013}}\label{alg:SVD}
  {\bf Input:} $A \in \C^{M\times N}$, $b\in\C^M$, $\eps > 0$ \\
  {\bf Output:} $x\in\C^N$ such that $Ax \approx b$
  \begin{algorithmic}[1]
    \State Compute the SVD, $A = U\Sigma V^*$ where
    \begin{equation*}
    U\Sigma V^* = \left(\begin{array}{cc}U_1 & U_2 \end{array} \right)     \left(\begin{array}{cc}\Sigma_1 & \\ & \Sigma_2 \end{array} \right) \left(\begin{array}{cc}V_1 & V_2 \end{array} \right)^*,
    \end{equation*}
      with $0 \leq \Sigma_2 < \eps I \leq \Sigma_1$.
    \State $x \gets V_1 \Sigma_1^{-1} U_1^* b$
  \end{algorithmic}
  \end{algorithm}

Algorithm~\ref{alg:SVD} is a standard method of solving a linear system using a truncated SVD \cite{Golub2013}. It is based on computing the full SVD of the matrix $A$, and subsequently discarding the singular values smaller than a threshold $\eps$. Assuming that $M = {\mathcal O}(N)$, this algorithm has cubic complexity ${\mathcal O}(N^3)$ even when the numerical rank of $A$ is small. For systems of low rank, the randomized algorithms that follow have a more favourable complexity. Nevertheless, because this algorithm is suggested in \cite{adcock2019frames,adcock2018fna2}, we prove bounds on the residual which are similar in flavour to those in \cite{adcock2019frames,adcock2018fna2}.

\begin{lemma}
  Let $x$ be computed by Algorithm \ref{alg:SVD}. Then
  \begin{equation*}
  \|b - Ax\|_2 \leq \inf_{v \in \C^N} \bigg\{\|b - Av \|_2 + \eps \cdot \|v\|_2 \bigg\}.
  \end{equation*}
  \begin{proof}
     We substitute $x = V_1 \Sigma_1^{-1} U_1^* b$ into the residual to obtain
     \begin{equation*}
      b - Ax = (I - AV_1 \Sigma_1^{-1} U_1^*)b.
     \end{equation*}
     One can expand the block form of the SVD of $A$ into $A = U_1\Sigma_1 V_1^* + U_2\Sigma_2 V_2^*$ and since the columns of $V$ are orthonormal vectors, we know that $V_2^* V_1 = 0$. Therefore,
     \begin{eqnarray*}
     AV_1 \Sigma_1^{-1} U_1^* &=& U_1\Sigma_1 V_1^*V_1 \Sigma_1^{-1} U_1^* + U_2\Sigma_2 V_2^*V_1 \Sigma_1^{-1} U_1^* \\
     &=& U_1\Sigma_1 \Sigma_1^{-1} U_1^* \\
     &=& U_1U_1^*.
     \end{eqnarray*}
    For any $v \in \mathbb{C}^{N}$, we can add and subtract $(I-U_1U_1^*)Av$ to get,
   \begin{eqnarray*}
     b - Ax &=& (I-U_1U_1^*) (b- Av) + (I-U_1U_1^*) A v.
     \end{eqnarray*}
   Since the columns of $U$ are orthonormal, we have $U_1^* U_2 = 0$ and $U_1^*U_1 = I$. Therefore,
   \begin{eqnarray*}
   b - Ax &=& (I-U_1U_1^*)(b-Av) + (I-U_1U_1^*)(U_1\Sigma_1 V_1^* + U_2 \Sigma_2 V_2^*)v \\
          &=& (I-U_1U_1^*)(b-Av) + U_2 \Sigma_2 V_2^* v.
   \end{eqnarray*}
   Since $\Sigma_2 < \eps I$ by assumption, $\|\Sigma_2\|_2 < \eps$. Also, $I-U_1U_1^*$ is an orthogonal projection onto the orthogonal complement of the range of $U_1$, so $\|I-U_1U_1^*\|_2 = 1$. Furthermore, $\|U_2\|_2  = \|V_2^*\|_2 = 1$, so the bound on the norm of the residual readily follows.
    \end{proof}
  \end{lemma}

Randomized algorithms based on matrix-vector products with random vectors can compute a truncated SVD at a lower cost in case the effective numerical rank is small. Algorithm~\ref{alg:randSVD} assumes that the user specifies a value $R=r+p$, where $r$ is the numerical rank of $A$ and larger values of $p$ yield higher probability of accuracy, as quantified in analysis that follows.

The existing analysis of randomized linear algebra (see \cite{halko2011finding} and references therein) focuses mostly on the accuracy of the matrix factorization. In this paper we are interested mainly in bounds of the residual of Algorithm~\ref{alg:randSVD}. To our knowledge, it is not possible to derive residual bounds directly from the error of the factorisation, so we derive the result directly using similar techniques to those in \cite{halko2011finding} for bounding expectation and tail probabilities of the random matrices involved.

\begin{algorithm}[h!]
  \caption{Randomized truncated SVD solver}\label{alg:randSVD}
  {\bf Input:} $A \in \C^{M\times N}$, $b\in\C^M$, $R \in \{1,\ldots,N\}$, $\eps > 0$ \\
  {\bf Output:} $x\in\C^N$ such that $Ax \approx b$
  \begin{algorithmic}[1]
    \State Generate $\Omega \sim \mathcal{N}(0,1;\R^{N\times R})$
    \State $\tilde{A} \gets A\Omega \in \C^{M\times R}$
    \State Compute the SVD, $\tilde A = \tilde U \tilde\Sigma \tilde V^*$ where
    \begin{equation*}
    \tilde U\tilde\Sigma \tilde V^* = \left(\begin{array}{cc}\tilde U_1 & \tilde U_2 \end{array} \right)     \left(\begin{array}{cc}\tilde \Sigma_1 & \\ & \tilde\Sigma_2 \end{array} \right) \left(\begin{array}{cc}\tilde V_1 & \tilde V_2 \end{array} \right)^*,
    \end{equation*}
    with $0 \leq \tilde\Sigma_2 < \eps I \leq \tilde\Sigma_1$. Here $\tilde U \in \C^{M\times R}, \tilde\Sigma \in \R^{R\times R}, \tilde V \in \C^{R\times R}$ but the dimensions of the blocks depend on the singular values.
    \State $y \gets \tilde V_1 \tilde\Sigma_1^{-1} \tilde U_1^* b$
    \State $x \gets \Omega y$
  \end{algorithmic}
\end{algorithm}

\begin{theorem}[Residual bounds for randomized truncated SVD solver]\label{thm:randomSVD}
  Assume that $A$ is such that $\rank_\eps(A) = r$ (as defined in equation \eqref{eqn:lowrank}) and let $x \in \C^{N}$ come from step 5 of Algorithm \ref{alg:randSVD} with $R = r + p$ for $p \geq 2$. Then
  \begin{equation*}
  \|b - Ax\|_2 \leq \inf_{v \in \C^N} \bigg\{\|b - Av \|_2 + \eps \cdot (1+\kappa) \cdot \|v\|_2 \bigg\},
  \end{equation*}
  where $\kappa$ is a non-negative-valued random variable satisfying
  \begin{equation*}
  \mathbb{E}\left\{ \kappa \right\} \leq 2\sqrt{\frac{r}{p-1}}, \qquad \mathbb{P}\left\{ \kappa > (2+u)\cdot s\cdot \sqrt{\frac{3r}{p+1}} \right\} \leq s^{-p} + \ee^{-\frac{u^2}{2}},
  \end{equation*}
for any $s \geq 1$, $u\geq 0$. Loosely speaking, $\kappa = \mathcal{O}\left( \sqrt{r} \right)$ with a high probability which improves rapidly as $p$ increases.
  \end{theorem}
\begin{remark}
  The probability distribution of $\kappa$ is similar to that of random variables appearing in the factorisation errors of \cite[Sec.~10]{halko2011finding}. Careful choices of $s$ and $u$ will give different bounds on the probability which depends on $p$. Following the example choices of parameters in \cite{halko2011finding}, setting $s = e$, $u = 2 + \sqrt{2p}$ and $p = 20$ in this Theorem shows that,
  \begin{equation*}
   \mathbb{E}\left\{ \kappa \right\} \leq 0.459 \sqrt{r}, \qquad \mathbb{P}\left\{ \kappa > 8.56\sqrt{r} \right\} \leq 4.13 \times 10^{-9}.
  \end{equation*}
  It might be tempting to let $p$ grow linearly with respect to $r$, since then the expected value of $\kappa$ is $\mathcal{O}(1)$ instead of $\mathcal{O}(\sqrt{r})$, but this is overkill, and we find that in practice a fixed $p$ such as $20$ works well.
  \end{remark}
  \begin{proof}
    Note that we can write $\tilde A = \tilde U_1 \tilde \Sigma_1 \tilde V_1^* + \tilde U_2 \tilde \Sigma_2 \tilde V_2^*$ with $\|\tilde \Sigma_2\|_2 \leq \eps$, by step 3 of Algorithm \ref{alg:SVD}. By assumption, $A$ has epsilon rank $r$, so its SVD is of the form
    \begin{equation*}
    A = \underbrace{U_1\Sigma_1V_1^*}_L + \underbrace{U_2 \Sigma_2 V_2^*}_E,
    \end{equation*}
    where $\Sigma_1 \in \R^{r\times r}$, the concatenated matrices $[U_1|U_2]$ and $[V_1|V_2]$ have orthonormal columns, and $\|E\|_\FF = \|\Sigma_2\|_\FF \leq \eps$. Substituting $x = \Omega\tilde V_1 \tilde\Sigma_1^{-1} \tilde U_1^* b$ into the residual, we get
    \begin{eqnarray*}
   b-Ax &=& b - \tilde{A}\tilde V_1 \tilde\Sigma_1^{-1} \tilde U_1^* b \nonumber \\
    &=& (I-\tilde U_1\tilde\Sigma_1 \tilde V_1^* \tilde V_1 \tilde\Sigma_1^{-1}\tilde U_1^* - \tilde U_2\tilde\Sigma_2 \tilde V_2^* \tilde V_1 \tilde\Sigma_1^{-1}\tilde U_1^*)b \\
    &=& (I-\tilde U_1 \tilde U_1^*)b,\nonumber
    \end{eqnarray*}
    by the identities for $\tilde{U}_1,\tilde{U}_2,\tilde{V}_1,\tilde{V}_2$ which follow from the orthonormal columns of $\tilde U$ and $\tilde V$. Now, for any $v \in \C^N$, if we write $b = (b-Av) + Av = (b-Av) + Ev + Lv$, then
    \begin{equation*}\label{eqn:bminusAx1}
   b-Ax = (I-\tilde U_1 \tilde U_1^*)(b-Av +Ev + Lv).
    \end{equation*}
    
    Consider the matrices $\Omega_1 = V_1^*\Omega$ and $\Omega_2 = V_2^*\Omega$. They are submatrices of the Gaussian matrix $[V_1 | V_2]^* \Omega$ (noting that the independence and Gaussian distribution of the elements is preserved by a unitary transformation). It follows that $\Omega_1$ and $\Omega_2$ are independent Gaussian matrices in $\R^{r\times (r+p)}$ and $\R^{(N-r)\times(r+p)}$ respectively. With probability 1, the rows of $\Omega_1$ are linearly independent, so that the pseudoinverse is in fact a right inverse i.e.~$\Omega_1 \Omega_1^\dagger = I_{r \times r}$ with probability 1. Therefore, we can write $L = U_1\Sigma_1 V_1^* = U_1\Sigma_1 \Omega_1\Omega_1^\dagger V_1^* = U_1\Sigma_1 V_1\Omega\Omega_1^\dagger V_1^* =  L \Omega\Omega_1^\dagger V_1^*$. Since $\tilde A = A\Omega$ (by step 2 of Algorithm \ref{alg:SVD}), this implies $L = (A-E)\Omega\Omega_1^\dagger V_1^* = \tilde A\Omega_1^\dagger V_1^* - U_2\Sigma_2\Omega_2\Omega_1^\dagger V_1^*$. Substituting this into equation \eqref{eqn:bminusAx1} gives,
    \begin{eqnarray}\label{eqn:bminusAx2}
    b-Ax &=& (I-\tilde U_1 \tilde U_1^*)(b-Av) \nonumber \\
    & & + (I-\tilde U_1 \tilde U_1^*)(E -U_2\Sigma_2\Omega_2\Omega_1^\dagger V_1^*)v \\
    & & + (I-\tilde U_1 \tilde U_1^*)\tilde{A}\Omega_1^\dagger V_1^*v \nonumber.
    \end{eqnarray}
    There are three terms here. Note that $\|I- \tilde U_1 \tilde U_1^*\|_2 \leq 1$ since $\tilde U_1$ has orthonormal columns, so the first term is bounded above by $\|b-Av\|_2$ and the second term is bounded above by $(\|E\|_{\FF} + \|\Sigma_2 \Omega_2 \Omega_1^\dagger\|_{\FF})\|v\|_2$. The third term requires more manipulation, as follows. Using the decomposition $\tilde{A} = \tilde U_1 \tilde \Sigma_1 \tilde V_1^* + \tilde U_2 \tilde \Sigma_2 \tilde V_2^*$ made at the start of the proof, we obtain the identity $(I-\tilde U_1 \tilde U_1^*)\tilde{A} =  \tilde U_2 \tilde \Sigma_2 \tilde V_2^*$. Using this, the third term in equation \eqref{eqn:bminusAx2} has norm that is readily confirmed to be bounded above by $\|\tilde\Sigma_2\|_2 \|\Omega_1^\dagger\|_\FF \|v\|_2$. 
    
    Combining our estimates for the three terms in equation \eqref{eqn:bminusAx2} provides the deterministic bound,
    \begin{eqnarray*}
    \|b - Ax \|_2 &\leq& \|b-Av\|_2 + \left(\|E\|_{\FF} + \|\Sigma_2\Omega_2\Omega_1^\dagger\|_\FF\right)\|v\|_2 + \|\tilde\Sigma_2\|_2\|\tilde \Omega_1^\dagger\|_\FF \|v\|_2 \\
    &\leq& \|b-Av\|_2 + \eps\left(1 + \|\eps^{-1}\Sigma_2\Omega_2\Omega_1^\dagger\|_\FF + \|\Omega_1^\dagger\|_\FF \right)\|v\|_2
    \end{eqnarray*}
    Now we bound the expectation and tail probabilities of the random variable
    \begin{equation*}
     \kappa = \|\eps^{-1}\Sigma_2\Omega_2\Omega_1^\dagger\|_\FF + \|\Omega_1^\dagger\|_\FF. 
     \end{equation*}
    Note that this has non-negative value, as required. From proposition \ref{prop:Gaussianbounds1} with $S^* = \eps^{-1}\Sigma_2$ and $T = \Omega_1^\dagger$ which is independent of the Gaussian matrix $\Omega_2$, we obtain
    \begin{equation*}\label{eqn:expectationbound}
    \mathbb{E}\left\{\|\eps^{-1}\Sigma_2\Omega_2\Omega_1^\dagger\|_\FF \quad | \quad \Omega_1 \right\} \leq \|\Omega_1^\dagger\|_\FF,
    \end{equation*}
    since $ \|\eps^{-1}\Sigma_2\|_\FF \leq 1$. Therefore $\mathbb{E}\left\{\kappa \right\} \leq 2 \mathbb{E}\left\{\|\Omega_1^\dagger\|_\FF \right\}$. Applying Proposition \ref{prop:Gaussianbounds2} to this yields the bound on the expectation of $\kappa$.
    
    Now we turn to the bounds on tail probabilities of $\kappa$. Following \cite[Thm.~10.8]{halko2011finding} with a crude simplification (in which we only consider Frobenius norms throughout the calculation), we condition on the event
    \begin{equation*}
    \mathcal{E}_s = \left\{ \Omega_1 : \|\Omega_1^\dagger\|_\FF < s \cdot \sqrt{\frac{3r}{p+1}} \right\},
    \end{equation*}
    where $s \geq 1$. Proposition \ref{prop:Gaussianbounds2} implies that the probability that $\mathcal{E}_s$ does not occur is $s^{-p}$. Conditional on $\mathcal{E}_s$, Proposition \ref{prop:Gaussianbounds1} gives us the inequality
    \begin{equation*}
    \mathbb{P}\left\{ \|\eps^{-1}\Sigma_2 \Omega_2 \Omega_1^\dagger\|_\FF > (1+u)\cdot s\cdot\sqrt{\frac{3r}{p+1}} \quad | \quad \mathcal{E}_s \right\} \leq \ee^{-\frac{u^2}{2}},
    \end{equation*}
    for all $u \geq 0$. This implies
    \begin{equation*}
    \mathbb{P}\left\{ \kappa > (2+u) \cdot s \cdot \sqrt{\frac{3r}{p+1}} \quad | \quad \mathcal{E}_s \right\} \leq \ee^{-\frac{u^2}{2}}.
    \end{equation*}
    Adding in the probability that $\mathcal{E}_s$ does not occur, in order to remove the conditioning, we arrive at the tail probability bound for $\kappa$.
    \end{proof}

\subsection{Truncated pivoted QR solvers}

Similar to the case of the SVD, first we formulate a standard algorithm based on a full pivoted QR decomposition in Algorithm~\ref{alg:QR}. A randomized algorithm with better complexity for systems of small numerical rank is Algorithm~\ref{alg:randQR}.

\begin{algorithm}[h!]
  \caption{Truncated pivoted QR solver \cite{Golub2013}}\label{alg:QR}
  {\bf Input:} $A \in \C^{M\times N}$, $b\in\C^M$, $r \in \{1,\ldots,N\}$ \\
  {\bf Output:} $x\in\C^N$ such that $Ax \approx b$
  \begin{algorithmic}[1]
    \State Compute a column pivoted QR decomposition $A \Pi = Q R$, with block forms,
    \begin{equation*}
    \Pi =  \left(\begin{array}{cc} {\Pi}_1 & {\Pi}_2 \end{array} \right), \quad {Q} = \left(\begin{array}{cc} {Q}_1 & {Q}_2 \end{array} \right), \quad {R} = \left(\begin{array}{cc} {R}_{11} & {R}_{12} \\ 0 & {R}_{22} \end{array} \right),
    \end{equation*}
    where $\Pi_1\in \C^{N\times r}$, $\Pi_2\in \C^{N\times N-r}$, $Q_1 \in \C^{M\times r}$, $Q_2 \in \C^{M\times {N-r}}$ $R_{11} \in \C^{r\times r}$, $R_{22} \in \C^{N- r \times N-r}$.
    \State $x \gets  \Pi_1 {R}_{11}^{-1}  Q_1^* b$
  \end{algorithmic}
\end{algorithm}

\begin{lemma}\label{lem:pivotedQR}
  Let $x$ be computed by Algorithm \ref{alg:QR}. Then
  \begin{equation*}
  \|b - A x\|_2 \leq \inf_{v \in \C^N} \bigg\{\|b - Av \|_2 + \|R_{22}\|_2 \cdot \|v\|_2 \bigg\}
  \end{equation*}
  \begin{proof}
    We substitute $x = \Pi_1 {R}_{11}^{-1}  Q_1^* b$ into the residual to obtain
    \begin{equation*}
    b - Ax = (I - A\Pi_1 {R}_{11}^{-1}  Q_1^*)b.
    \end{equation*}
    Note that $A\Pi_1 = Q_1 R_{11}$. Therefore, $A\Pi_1 {R}_{11}^{-1}  Q_1^* = Q_1 Q_1^*$ and so $b-Ax = (I-Q_1Q_1^*)b$. For any $v\in\C^N$ we add and subtract $(I-Q_1Q_1^*)Av$ to obtain
    \begin{eqnarray*}
      b - Ax &=& (I-Q_1Q_1^*)(b- Av) + (I-Q_1Q_1^*) A v.
    \end{eqnarray*}
 Let us deal with the $(I-Q_1Q_1^*) A v$ term. Note that $\Pi$ is merely a permutation matrix, so $\Pi \Pi^{\mathrm{T}} = I$. Therefore,
 \begin{eqnarray*} 
   (I-Q_1Q_1^*) A v &=& (I-Q_1Q_1^*) A\Pi \Pi^{\mathrm{T}}v \\
                    &=& (I-Q_1Q_1^*) QR \Pi^{\mathrm{T}}v \\
                    &=& \left(\begin{array}{cc} 0 & {Q}_2 \end{array} \right) R\Pi^{\mathrm{T}}v \\
                    &=& \left(\begin{array}{cc} 0 & {Q}_2R_{22} \end{array} \right) \Pi^{\mathrm{T}}v \\
                    &=& {Q}_2R_{22} \Pi_2^{\mathrm{T}} v.
   \end{eqnarray*}
  Therefore, $b-Ax = (I-Q_1Q_1^*)(b- Av) + {Q}_2R_{22} \Pi_2^{\mathrm{T}} v$, and the bound follows readily from $\|I-Q_1Q_1^*\|_2 \leq 1$, $\|Q_2\|_2\leq 1$ and $\|\Pi_2^{\mathrm{T}}\|_2 \leq 1$.
    \end{proof}
  \end{lemma}

\begin{algorithm}[h!]
  \caption{Randomized truncated pivoted QR solver}\label{alg:randQR}
  {\bf Input:} $A \in \C^{M\times N}$, $b\in\C^M$, $R \in \{1,\ldots,N\}$, $\eps > 0$ \\
  {\bf Output:} $x\in\C^N$ such that $Ax \approx b$
  \begin{algorithmic}[1]
    \State Generate $\Omega \sim \mathcal{N}(0,1;\R^{N\times R})$
    \State $\tilde{A} \gets A\Omega \in \C^{M\times R}$
    \State Compute a column pivoted QR decomposition $\tilde{A}\tilde{\Pi} = \tilde{Q}\tilde{R}$, with block forms,
    \begin{equation*}
     \tilde\Pi =  \left(\begin{array}{cc} {\tilde\Pi}_1 & {\tilde\Pi}_2 \end{array} \right), \quad \tilde{Q} = \left(\begin{array}{cc} \tilde{Q}_1 & \tilde{Q}_2 \end{array} \right), \quad \tilde{R} = \left(\begin{array}{cc} \tilde{R}_{11} & \tilde{R}_{12} \\ 0 & \tilde{R}_{22} \end{array} \right),
    \end{equation*}
  with $0 \leq \diag(\tilde R_{22}) < \eps I \leq \diag( \tilde R_{11})$. Here $\tilde \Pi \in \R^{M\times R}, \tilde Q \in \C^{M\times R}, \tilde R \in \C^{R\times R}$ but the dimensions of the blocks depend on the diagonal entries of $\tilde{R}$.
    \State $y \gets \tilde \Pi_1 \tilde{R}_{11}^{-1} \tilde Q_1^* b$
    \State $x \gets \Omega y$
  \end{algorithmic}
\end{algorithm}

\begin{theorem}[Residual bounds for randomized pivoted QR solver]\label{thm:randomQR}
  Assume that $A$ is such that $\rank_\eps(A) = r$ (as defined in equation \eqref{eqn:lowrank}) and let $x \in \C^{N}$ come from step 5 of Algorithm \ref{alg:randQR} with $R = r + p$ for $p \geq 2$. Then
  \begin{equation*}
  \|b - Ax\|_2 \leq \inf_{v \in \C^N} \bigg\{\|b - Av \|_2 + \eps \cdot (1+\kappa) \cdot \|v\|_2 \bigg\},
  \end{equation*}
  where $\kappa$ is a non-negative-valued random variable satisfying
  \begin{equation*}
  \mathbb{E}\left\{ \kappa \right\} \leq (1+\sqrt{r+p})\sqrt{\frac{r}{p-1}}, \qquad \mathbb{P}\left\{ \kappa > (1 + \sqrt{r+p} + u)\cdot s\cdot \sqrt{\frac{3r}{p+1}} \right\} \leq s^{-p} + \ee^{-\frac{u^2}{2}},
  \end{equation*}
  for any $s \geq 1$, $u\geq 0$. Loosely speaking, $\kappa = \mathcal{O}\left( r \right)$ with a high probability which improves rapidly as $p$ increases.
\end{theorem}
\begin{remark}
  The bounds here on $\kappa$ are larger than those found for the randomized truncated SVD algorithm in Theorem \ref{thm:randomSVD}. This is related to the possibility that the diagonals of $\tilde{R}_{22}$ in the pivoted QR algorithm can bear no relation to the underlying numerical rank of $\tilde{A}$. This is discussed in \cite[Sec.~5.4.3]{Golub2013}, and therein it is noted that: ``Nevertheless, in practice, small trailing $R$-submatrices almost always emerge that correlate well with the underlying [numerical] rank''. Hence, in practice, we would expect $\kappa = \mathcal{O}(\sqrt{r})$ as in the randomized truncated SVD algorithm.
  \end{remark}
\begin{proof}
 By assumption, $A$ has epsilon rank $r$, so its SVD is of the form
  \begin{equation*}
  A = \underbrace{U_1\Sigma_1V_1^*}_L + \underbrace{U_2 \Sigma_2 V_2^*}_E,
  \end{equation*}
  where $\Sigma_1 \in \R^{r\times r}$, the concatenated matrices $[U_1|U_2]$ and $[V_1|V_2]$ have orthonormal columns, and $\|E\|_\FF = \|\Sigma_2\|_\FF \leq \eps$. Substituting $x = \Omega\tilde \Pi_1 \tilde R_{11}^{-1} \tilde Q_1^* b$ into the residual gives
  \begin{eqnarray*}
  b-Ax &=& b - \tilde{A}\tilde \Pi_1 \tilde R_{11}^{-1} \tilde Q_1^* b.
  \end{eqnarray*}
  Note that $\tilde A\tilde\Pi_1 = \tilde Q_1 \tilde R_{11}$, therefore, $\tilde A\tilde \Pi_1 {\tilde R}_{11}^{-1}  \tilde Q_1^* = \tilde Q_1 \tilde Q_1^*$, so $b - Ax = (I- \tilde{Q}_1\tilde{Q}_1^*)b$.
 Now, for any $v \in \C^N$, if we write $b = (b-Av) + Av = (b-Av) + Ev + Lv$, then
  \begin{equation*}\label{eqn:bminusAx1QR}
  b-Ax = (I-\tilde Q_1 \tilde Q_1^*)(b-Av +Ev + Lv).
  \end{equation*}
  
  Using the same notation and result derived in the proof of Theorem \ref{thm:randomSVD}, we have $L = \tilde A\Omega_1^\dagger V_1^* - U_2\Sigma_2\Omega_2\Omega_1^\dagger V_1^*$. Substituting this into equation~\eqref{eqn:bminusAx1QR} gives,
  \begin{eqnarray*}\label{eqn:bminusAx2QR}
  b-Ax &=& (I-\tilde Q_1 \tilde Q_1^*)(b-Av) \nonumber \\
  & & + (I-\tilde Q_1 \tilde Q_1^*)(E -U_2\Sigma_2\Omega_2\Omega_1^\dagger V_1^*)v \\
  & & + (I-\tilde Q_1 \tilde Q_1^*)\tilde{A}\Omega_1^\dagger V_1^*v \nonumber.
  \end{eqnarray*}
  There are three terms here of which the 2-norm needs to be bounded. Note that $\|I- \tilde Q_1 \tilde Q_1^*\|_2 \leq 1$ since $\tilde Q_1$ has orthonormal columns, so the first term is bounded above by $\|b-Av\|_2$ and the second term is bounded above by $(\|E\|_{\FF} + \|\Sigma_2 \Omega_2 \Omega_1^\dagger\|_{\FF})\|v\|_2$. The third term requires more manipulation, as follows. Following the same reasoning as in the proof of Lemma \ref{lem:pivotedQR}, we have $(I-\tilde Q_1 \tilde Q_1^*)\tilde{A} =  {\tilde Q}_2 \tilde R_{22} \tilde \Pi_2^{\mathrm{T}}$. Using this, the third term in equation \eqref{eqn:bminusAx2} has norm that is readily confirmed to be bounded above by $\|\tilde R_{22}\|_{\FF} \|\Omega_1^\dagger\|_\FF \|v\|_2$.
  
  Combining our estimates for the three terms in equation \eqref{eqn:bminusAx2} provides the deterministic bound,
  \begin{eqnarray*}
    \|b - Ax \|_2 &\leq& \|b-Av\|_2 + \left(\|E\|_{\FF} + \|\Sigma_2\Omega_2\Omega_1^\dagger\|_\FF\right)\|v\|_2 + \|\tilde R_{22}\|_{\FF}\|\tilde \Omega_1^\dagger\|_\FF \|v\|_2 \\
    &\leq& \|b-Av\|_2 + \eps\left(1 + \|\eps^{-1}\Sigma_2\Omega_2\Omega_1^\dagger\|_\FF + \|
    \eps^{-1} \tilde R_{22}\|_{\FF} \|\Omega_1^\dagger\|_\FF \right)\|v\|_2
  \end{eqnarray*}
 
We can bound $\|\eps^{-1} \tilde R_{22}\|_{\FF}$ by $\sqrt{r+p}$ by the following argument. Column pivoting ensures that the first column of $\tilde{R}_{22}$ has the greatest norm of all its columns \cite[Sec.~5.4.2]{Golub2013}. Therefore, since $\tilde{R}_{22}$ has at most $r+p$ columns, we have
\begin{equation*}
  \|\tilde{R}_{22}\|_{\FF}^2 \leq (r+p) \left|\left[\tilde{R}_{22}\right]_{11}\right|^2 \leq (r+p)\eps^2.
  \end{equation*}

 The proof is completed if we can bound the expectation and tail probabilities of the random variable
\begin{equation*}
  \kappa = \|\eps^{-1}\Sigma_2\Omega_2\Omega_1^\dagger\|_\FF + \sqrt{r+p} \|\Omega_1^\dagger\|_\FF,
\end{equation*}
with the appropriate bounds. This can be done by almost exactly the same analysis as in the proof of Theorem \ref{thm:randomSVD}.
\end{proof}

\subsection{Computational cost of the randomized solvers}

Both Algorithm \ref{alg:randSVD} and Algorithm \ref{alg:randQR} consist of 5 steps. We break down the computational cost for each step. Recall that the input parameters are $A \in \C^{M\times N}$, $b\in \C^M$, $R \in \{1,\ldots, N\}$, $\eps > 0$. We will describe the costs using big $\mathcal{O}$ notation while regarding $A$, $b$ and $\eps$ as fixed quantities whose properties are hidden in the big $\mathcal{O}$. The number of operations required to apply $A$ to a vector is denoted $\mathrm{T}_{\mathrm{mult},A}$.

\begin{enumerate}
  \item Generate $N\cdot R$ Gaussian random numbers: $\mathcal{O}(N\cdot R)$
  \item Apply the matrix $A$ to $R$ vectors: $\mathcal{O}(R\cdot \mathrm{T}_{\mathrm{mult},A})$
  \item Compute the SVD or pivoted QR factorization of an $M \times R$ matrix: $\mathcal{O}(M\cdot R^2)$
  \item Apply a sequence of matrices whose with one dimension smaller than $R$ and the other smaller than $M$: $\mathcal{O}(R\cdot M)$
  \item Apply a matrix of size $N \times R$ to a vector: $\mathcal{O}(R\cdot N)$.
  \end{enumerate}
The total computational cost in floating point operations is therefore,
\begin{equation*}
\mathcal{O}(R \mathrm{T}_{\mathrm{mult},A} + R^2 M).
\end{equation*}
Hence, these algorithms have improved computational cost when compared to their non-randomized counterparts, Algorithm \ref{alg:SVD} and Algorithm \ref{alg:QR} (which require $\mathcal{O}(MN^2)$ operations), when $R = \mathrm{o}(N)$, with further gains to be had if $A$ supports fast matrix-vector multiplication.
% !TEX root = AZalgorithmV3.tex
\section{Fourier extension}
\label{s:fourext}

We proceed by illustrating the AZ algorithm for a number of cases. In particular, we show how to identify a suitable $Z$ matrix in different settings having to do with function approximation. The matrix $Z^*$ should be close to a pseudo-inverse of $A$ up to a numerically low rank error -- recall Lemma~\ref{lem:splitting}. For function approximation, in practice this means that $Z^*b$ are coefficients yielding an accurate approximation for a space of functions whose complement is of low dimension (relative to the degree of the approximation).

The first case is that of function approximation using Fourier extension frames. Here it is known that, up to a real-valued normalization factor, $Z=A$: see \cite[Algorithm 2]{matthysen2016fastfe} and \cite[Algorithm 1]{matthysen2017fastfe2d}. Compared to these references, we emphasize here how $Z$ arises from the so-called canonical dual frame, which yields extremely accurate results for functions compactly supported within the inner domain $\Omega$. The complement of this space is of low dimension related to the Hausdorff measure of the boundary \cite{matthysen2017fastfe2d}..

\subsection{Approximation with the Fourier extension frame}

We refer the rea\-der to \cite{christensen2008framesandbases} for an introduction to frames, and to \cite{adcock2019frames,adcock2018fna2} for an exposition on their use in function approximation. We merely recall that a system of functions $\Phi = \{ \phi_k\}_{k=1}^\infty$ is a frame for a separable Hilbert space $H$ if it satisfies the frame condition: there exist constants $0 < A,B < \infty$ such that
\begin{equation}\label{eq:framecondition}
  A \Vert f \Vert_H^2 \leq \sum_{k=1}^\infty \langle f, \phi_k \rangle_H^2 \leq B \Vert f \Vert_H^2.
\end{equation}
For comparison, note that for an orthonormal basis \eqref{eq:framecondition} holds with $A=B=1$. However, unlike a basis, a frame may be redundant.

\begin{figure}
 \begin{center}
  \subfigure[A punctured disk]{
   \includegraphics[height=4cm]{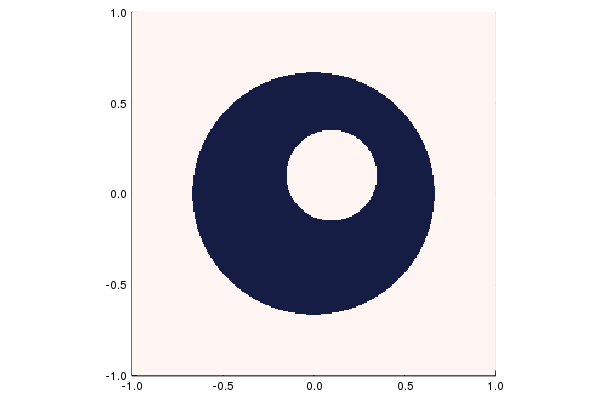}    }
    \hspace{0.5cm}
  \subfigure[Approximation of $\ee^{x+y}$]{
   \includegraphics[height=4cm]{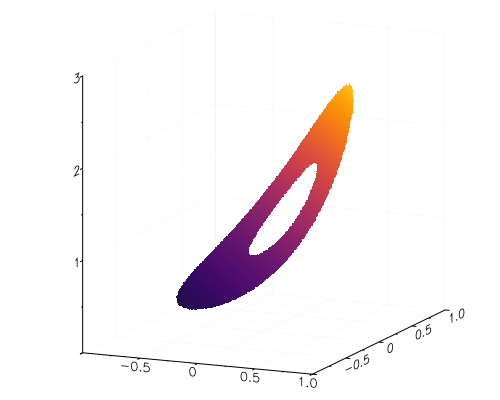}
    }
 \end{center}
 \caption{Fourier extension approximation of the function $f(x,y)=\ee^{x+y}$ on a punctured disk. The function is represented by a periodic Fourier series in the bounding box $[-1,1]^2$.}\label{fig:fourext_disk}
\end{figure}

To fix one frame, we consider a $d$-dimensional orthonormal Fourier basis $\Phi$ on the box $[-1,1]^d$,
\[
\Phi = \{ \phi_\mathbf{n} \}_{\mathbf{n} \in \Z^d}, \qquad \mbox{with} \qquad \phi_{\mathbf{n}}(\mathbf{x}) = 2^{-d/2} \ee^{\ii \pi \mathbf{n} \cdot \mathbf{x}}, \qquad \mathbf{x} \in [-1,1]^d.
\]
The restriction of this basis to a smaller domain $\Omega \subset [-1,1]^d$ yields a Fourier extension frame for $L^2(\Omega)$ (see \cite{huybrechs2010fourierextension} and \cite[Example 1]{adcock2019frames}). Any smooth function $f$ defined on $\Omega$ is well approximated in this frame, even if it is not periodic (which not even be a well-defined concept for an irregularly shaped $\Omega$). Indeed, each approximation in the frame implicitly corresponds to a \emph{periodic extension} of $f$ from $\Omega$ to $[-1,1]^d$. An example is shown in Fig.~\ref{fig:fourext_disk}.

The approximation problem requires the truncation of the infinite frame. Thus, assuming that $N^{1/d}$ is an integer, we define the subset of $N$ functions
\[
\Phi_N = \{ \phi_\mathbf{n}(\mathbf{x}) \}_{\mathbf{n} \in I_N},
\]
with
\[
  I_N = \left\{\mathbf{n} = (n_1,n_2,\ldots,n_d) \in \Z^d : -\frac12 N^{1/d} \leq n_1,\ldots,n_d \leq \frac12 N^{1/d} \right\},
\]
where we consider only $\mathbf{x} \in \Omega$. The best approximation to $f$ in the norm of $L^2(\Omega)$ leads to a linear system that involves the Gram matrix $G$ of $\Phi_N$ with elements
\begin{equation}\label{eq:gram_element}
G_{\mathbf{n},\mathbf{m}} = \langle \phi_{\mathbf{n}}, \phi_{\mathbf{m}} \rangle_{L^2(\Omega)}, \qquad \mathbf{n},\mathbf{m} \in I_N.
\end{equation}
Since inner products on $\Omega$ may be difficult to compute, especially for multidimensional domains, we resort to a \emph{discrete least squares approximation} instead. Given a set of points $\{ \mathbf{x}_m \}_{m=1}^M \subset \Omega$, this leads to the linear system $Ax=B$ with\footnote{The index $\mathbf{n}$ of $\phi_{\mathbf{n}}$ in the definition of the truncated frame $\Phi_N$ is a multi-index. For simplicity, in this expression we have identified each multi-index $\mathbf{n} \in I_N$ with an integer $n$ between $1$ and $N$.}
\begin{equation}\label{eq:fourext_leastsquares}
 A_{m,n} = \phi_n(\mathbf{x}_m), \qquad \mbox{and} \qquad B_m = f(\mathbf{x}_m), \qquad n=1,\ldots,N, \quad m=1,\ldots,M.
\end{equation}
It follows from the theory in \cite{adcock2019frames,adcock2018fna2} that this linear system has solutions that are stable least squares fits, and one such solution can be found using a truncated (and thereby regularized) singular value decomposition of $A$.

\subsection{The canonical dual frame}

The concept of dual frames does not play a large role in the analysis in \cite{adcock2019frames,adcock2018fna2}, but it is important in the theory of frames itself and, as we shall see, it is also relevant for the application of the AZ algorithm. A frame $\Phi = \{ \phi_k \}_{k=1}^\infty$ for a Hilbert space $H$ has a \emph{dual frame} $\Psi = \{ \psi_k \}_{k=1}^\infty$ if
\begin{equation}\label{eq:dualframe}
 f = \sum_{k=1}^\infty \langle f, \psi_k \rangle \phi_k, \qquad \forall f \in H.
\end{equation}
Convergence is understood in the norm of $H$. The above expansion simplifies if $\Phi$ is an orthonormal basis: in that case $\Psi = \Phi$. Yet, since frames may be redundant, it is important to remark that \eqref{eq:dualframe} may not be the only representation of $f$ in the frame. Correspondingly, dual frames are not necessarily unique.

The so-called \emph{canonical dual frame} plays a special role among all dual frames. For the case of Fourier extension, this dual frame corresponds to the Fourier series on $[-1,1]^d$ of $\overline{f}$, defined as the extension of $f$ by zero. The canonical dual frame expansion is readily obtained and easily truncated. For example, one can sample $\overline{f}$ in a regular grid on $[-1,1]^d$ and use the FFT. Unfortunately, though the expansion converges in norm, it does not actually yield an accurate approximation for most $f$. Since extension by zero introduces a discontinuity along the boundary of $\Omega$, in general the Fourier series of $\overline{f}$ exhibits the Gibbs phenomenon.

Nevertheless, the canonical dual frame expansion does lead to accurate approximations for a large class of functions, namely those that are compactly support on $\Omega$. In that case, their extension by zero results in a smooth and periodic function, and rapid convergence of the Fourier series of $\overline{f}$ follows. Thus, we have a simple and efficient FFT-based solver for a large subclass of functions, from which we can derive a suitable matrix $Z$.

\subsection{Characterizing the Z-matrix: continuous duality}

\begin{figure}
 \begin{center}
  \subfigure[Approximation of $\ee^x$ (solid blue) and its periodic Fourier extension (dashed red)]{
   \includegraphics[height=3.2cm]{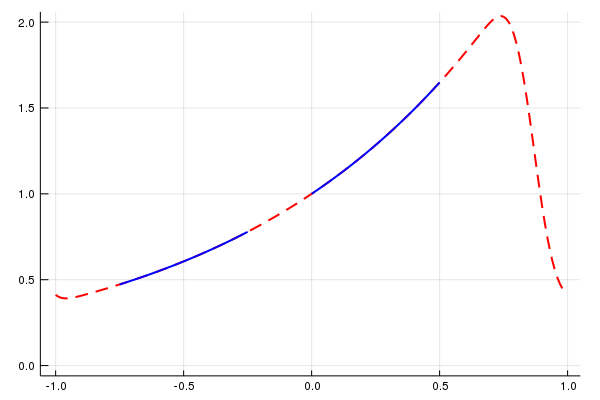}
    }
    \hspace{0.2cm}
  \subfigure[Singular values of the Gram matrix]{
   \includegraphics[height=3.2cm]{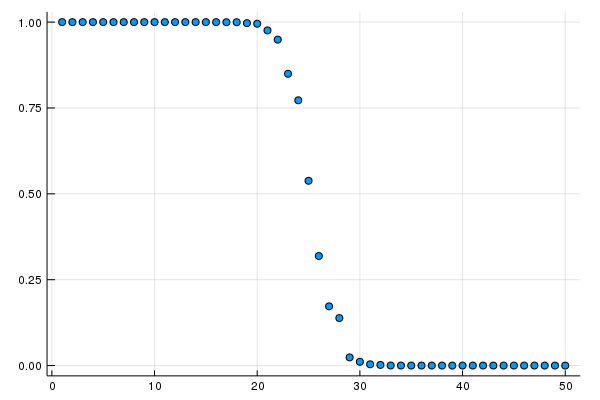}
    }
 \end{center}
 \caption{Fourier extension approximation of the function $f(x)=\ee^x$ on the domain $[-\frac34,-\frac14] \cup [0, \frac12]$, using orthonormal Fourier series on $[-1,1]$ ($N=50$). The best approximation was computed using the Gram matrix \eqref{eq:gram_element}. Its singular values cluster exponentially near $1$ and near $0$.}\label{fig:fourext_interval}
\end{figure}

We will first illustrate the properties of an approximation computed using the Gram matrix \eqref{eq:gram_element}. We consider a fully discrete approximation next. In Fig.~\ref{fig:fourext_interval} the computational domain is $\Omega = [-\frac34,-\frac14] \cup [0, \frac12]$. The frame is the restriction of the normalized Fourier basis on $[-1,1]$ to that domain. The left panel shows the approximation to $\ee^x$ on this non-connected domain. The right panel shows the singular values of the Gram matrix for $N=50$. The singular values cluster near $1$ and near $0$.\footnote{The properties of the Gram matrix are well understood for the case where $\Omega$ is a regular subinterval of $[-1,1]$. In that case, the Gram matrix is also known as the \emph{prolate matrix} \cite{slepian1978prolateV}. For example, it is known that the size of the cluster near $1$ is determined by the size of $\Omega$ relative to that of $[-1,1]$. Moreover, the clustering is exponential \cite{edelman1999futurefft}. The experiment shows that these properties are preserved for a non-connected subset $\Omega$.}

Any matrix $A$ with clustered singular values as shown in Fig.~\ref{fig:fourext_interval}(b) is amenable to the AZ algorithm with the simple choice $Z=A$, regardless of the underlying application. Indeed, if $A$ has left and right singular vectors $u_k$ and $v_k$ respectively, with singular value $\sigma_k$, one has that
\[
Z^* Av_k = V \Sigma U^* U \Sigma V^* v_k = \sigma_k^2 v_k.
\]
Thus, if $\sigma_k \approx 1$, we have that $Z^* Av_k \approx v_k$. Similarly, $A Z^* u_k \approx u_k$. %One can think of $Z$ as the (approximate) pseudo-inverse of $A$ on the subspace spanned by the first cluster of singular vectors, i.e., the singular vectors corresponding to singular values close to $1$.
Since the singular values of $A$ decay quickly away from the cluster near $1$, up to a small tolerance and low rank error, $Z$ acts like a generalized inverse of $A$.

The clustering of singular values precisely near the unit value $1$ is a consequence of the normalization of the Fourier series on $[-1,1]$. A different normalization of the basis functions is easily accommodated, since it only results in a diagonal scaling of the Gram matrix. In that case matrix $Z$ may be modified using the inverse of that scaling. That is, starting from a valid combination of $A$ and $Z$, we can use the combination
\[
A_1 = DA, \qquad Z_1 = (D^{-1})^*Z,
\]
so that their product $Z_1^* A_1 = Z^* D^{-1} D A = Z^* A$ remains unchanged.

\subsection{Characterizing the Z-matrix: discrete duality}

The discrete linear system~\eqref{eq:fourext_leastsquares} represents a more significant change to the linear system than a mere renormalization of the Gram matrix. The matrix elements are not given in terms of inner products but in terms of discrete function evaluations. Thus, we have to consider duality with respect to evaluation in the discrete grid. Fortunately, Fourier series are orthogonal with respect to evaluations on a periodic equispaced grid. Consider a set of $L$ normalized and univariate Fourier basis functions and their associated (periodic) equispaced grid $\{ x_l \}_{l=1}^L$, then we have:
\begin{equation}\label{eq:discrete_orthogonality}
 \sum_{l=1}^L \phi_i(x_l) \overline{\phi_j(x_l)} = L \, \delta_{i-j}, \qquad 1 \leq i,j \leq L.
\end{equation}
Equivalently, as is well known, the full DFT matrix $F$ of size $L \times L$ with entries $\phi_i(x_l)$ has inverse $F^*/L$. The discrete least squares problem \eqref{eq:fourext_leastsquares} follows by choosing a subset of $M < L$ points that belong to a subdomain $\Omega$ and a subset of $N < M$ basis functions. Thus, $A \in \mathbb{C}^{M \times N}$ is a submatrix of $F$, with rows corresponding to the selected points and columns corresponding to the selected basis functions. In this discrete setting, owing to \eqref{eq:discrete_orthogonality}, we can choose
\begin{equation}
Z = A / L.
\end{equation}

An interpretation can be given as follows. The multiplication with $A$ corresponds to an extension from $N$ to $L$ Fourier coefficients, followed by the DFT of length $L$, and followed by the restriction to $M$ points in $\Omega$. Multiplication with $Z^*$ corresponds to extension by zero-padding in the time-domain, followed by the inverse DFT, and restriction in the frequency domain. The vector $x = Z^* b$ is an accurate solution of the system $Ax=b$ in the special case where the sampled function is compactly supported on $\Omega$. Indeed, in that case zero-padding results in a smooth and periodic function, for which the discrete inverse Fourier transform gives an accurate approximation.
\begin{figure}
	\begin{center}
		\subfigure[Singular values of $A$]{
			\includegraphics[width=.4\textwidth]{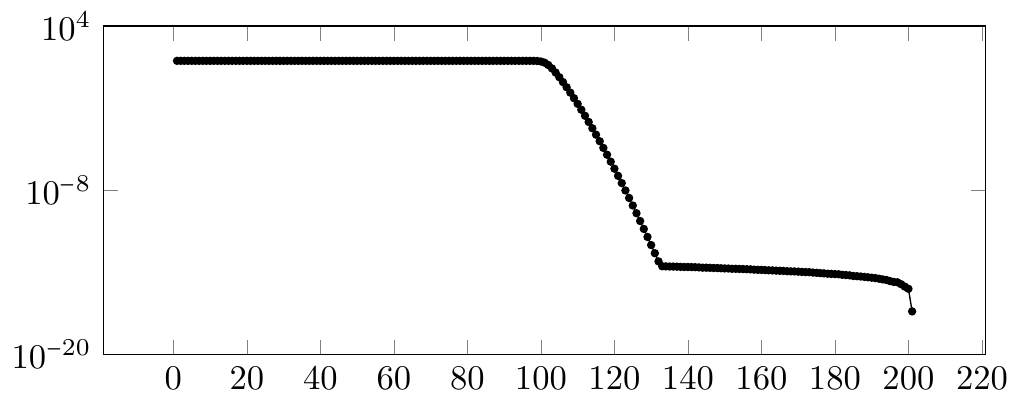}
		}
		\hspace{0.5cm}
		\subfigure[Singular values of $Z^*$]{
			\includegraphics[width=.4\textwidth]{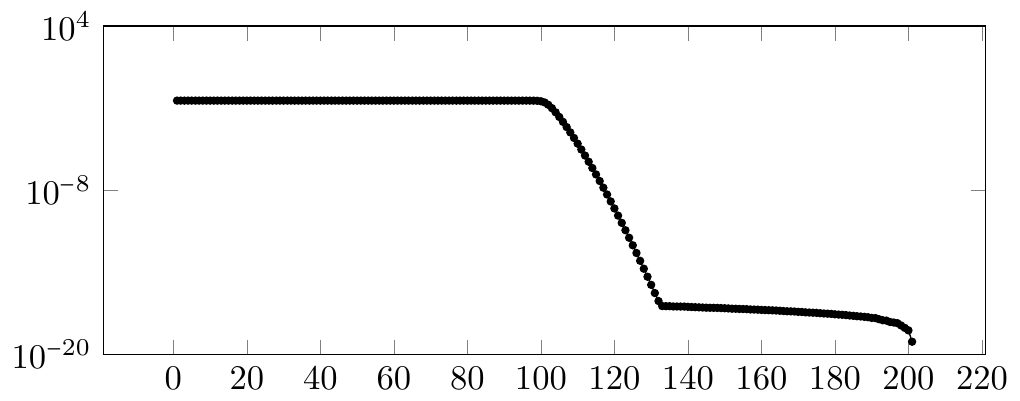}
		}\\
		\subfigure[Singular values of $A-AZ^*A$]{
			\includegraphics[width=.4\textwidth]{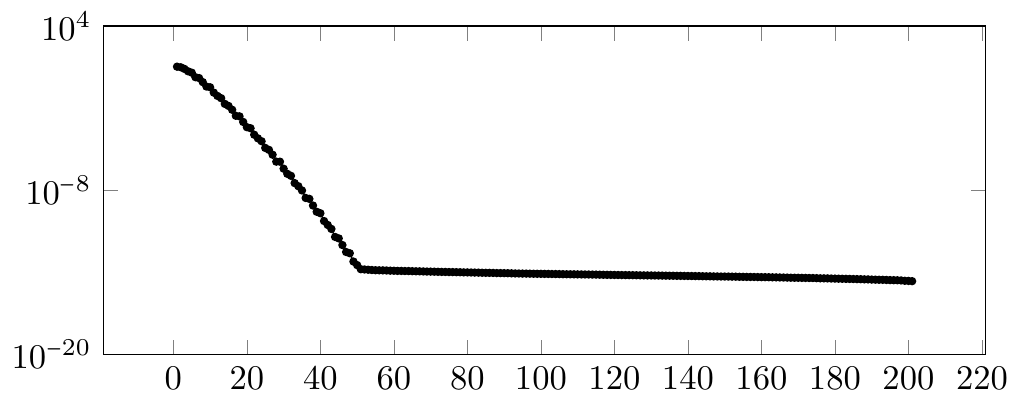}
		}
		\hspace{0.5cm}
		\subfigure[Timings in seconds]{
			\includegraphics[width=.4\textwidth]{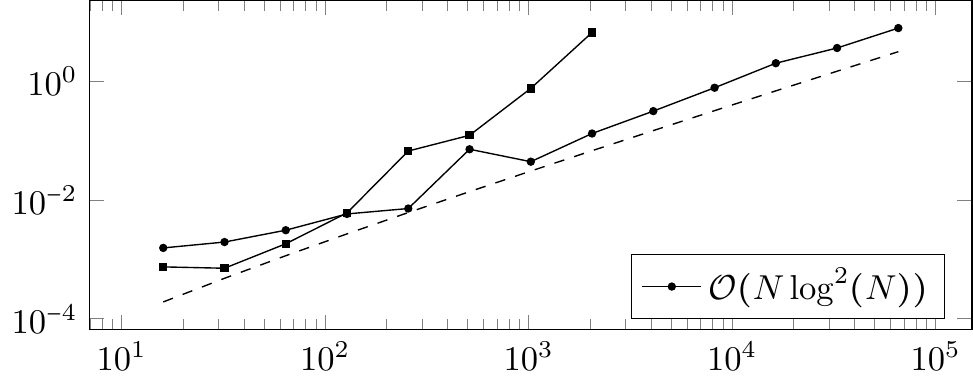}
		}
	\end{center}
	\caption{\label{fig:fourier_az}Discrete Fourier extension approximation  on the domain $\left[-\tfrac 12,\tfrac 12\right]$, using orthonormal Fourier series on $[-1,1]$. In the three first panels $N=201$. The last panel shows the timings in seconds for calculating the approximation to $f(x)=x$ for $N=16,32,\dots,2^{16}$ using the AZ-algorithm with a randomized SVD (dots), and a direct solver (squares). The dashed line shows $\mathcal O(N\log^2(N))$.}
\end{figure}

In Figure \ref{fig:fourier_az} we show the singular value patterns of $A$, $Z^*$, and $A-AZ^*A$ for the discrete Fourier extension from $\left[-\tfrac 12,\tfrac 12\right]$ to $[-1,1]$, with $N=201$, $M \approx 2N$ and $L \approx 2M$. The singular values of $A$ cluster near $\sqrt{L}$, and those of $Z=A/L$ near $1/\sqrt{L}$. The rank of the plunge region is known to be $\mathcal O(\log( N))$. The matrix $A-AZ^*A$ isolates this plunge region and thus has rank $\mathcal O(\log( N))$. Since applying $A$ and $Z^*$ is $\mathcal O(N\log( N))$ while the rank of the problem in the first step is $\mathcal O(\log( N))$, the overall computational complexity of the AZ algorithm is in this case $\mathcal O(N\log^2(N))$.

The setting is entirely the same in 2D, but the plunge region is larger relative to the size of the overall approximation problem \cite{matthysen2017fastfe2d}. The AZ algorithm was used to produce the 2D example in Fig.~\ref{fig:fourext_disk}.

\subsection{Generalization to other extension frames}

\begin{figure}
	\begin{center}
		\subfigure[Singular values of $A$]{
			\includegraphics[width=.4\textwidth]{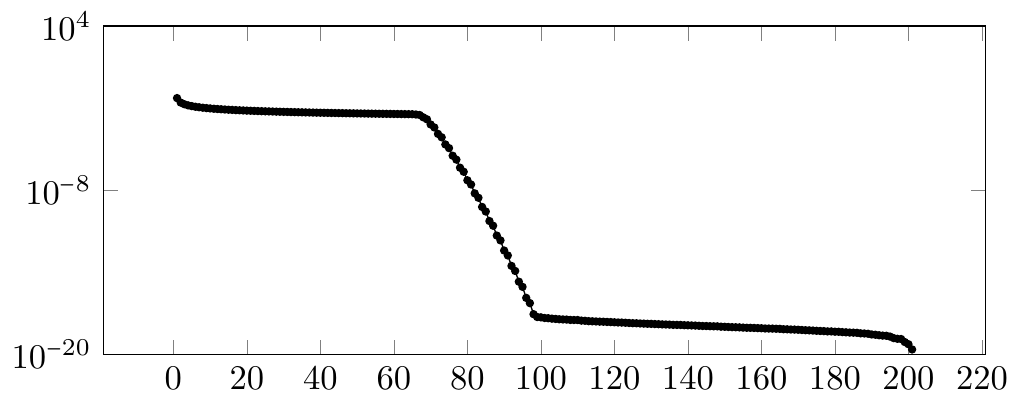}
		}
		\hspace{0.5cm}
		\subfigure[Singular values of $Z^*$]{
			\includegraphics[width=.4\textwidth]{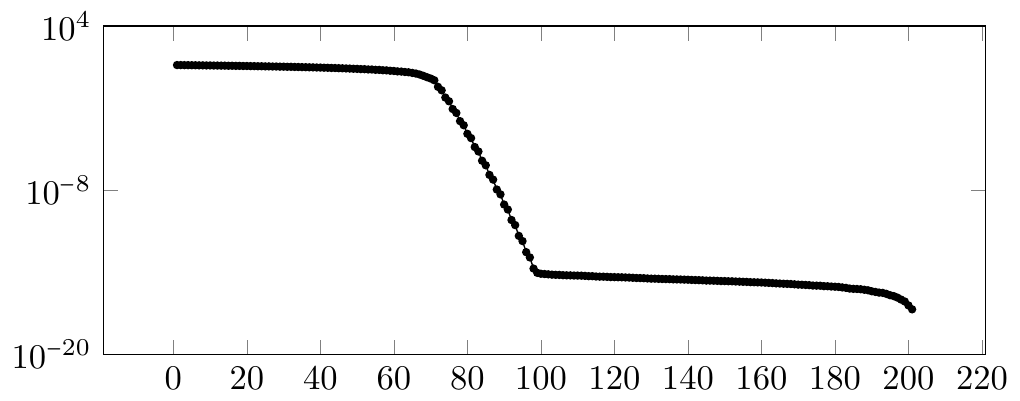}
		}\\
		\subfigure[Singular values of $A-AZ^*A$ (Legendre)]{
			\includegraphics[width=.4\textwidth]{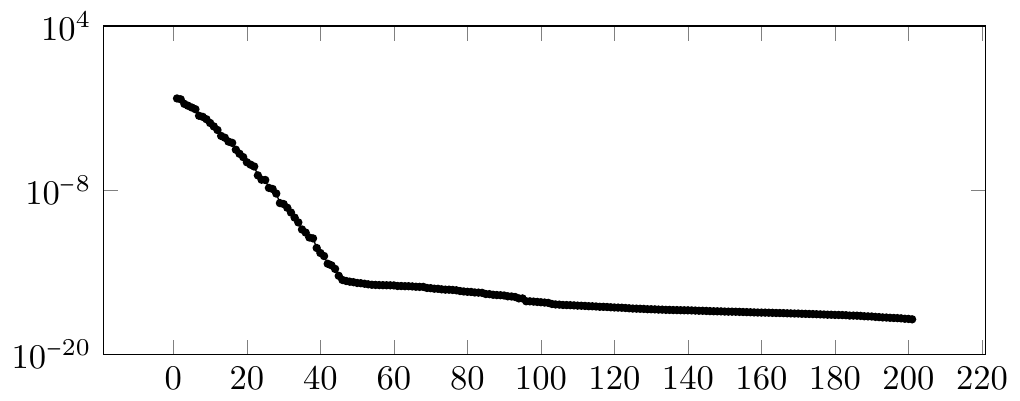}
		}
		\hspace{0.5cm}
		\subfigure[Singular values of $A-AZ^*A$ (Chebyshev)]{
			\includegraphics[width=.4\textwidth]{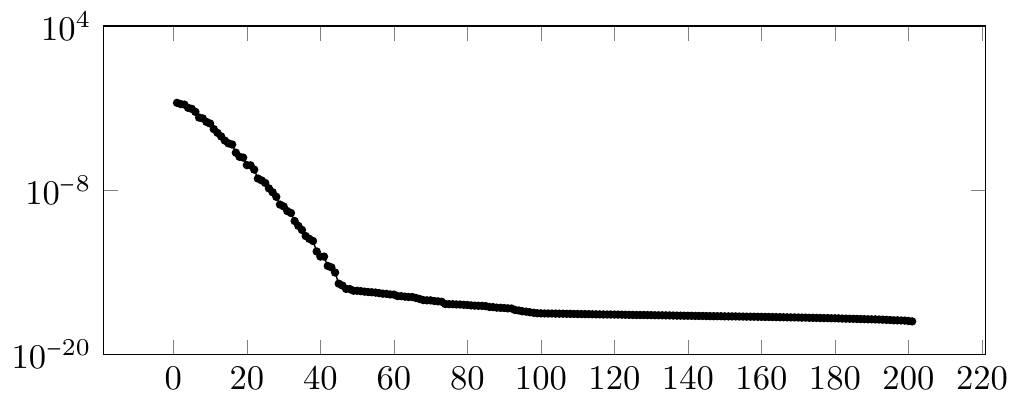}
		}
	\end{center}
	\caption{\label{fig:legendre_az}Legendre extension approximation on the domain $\left[-\tfrac12,\tfrac12\right]$, using classical Legendre series on $[-1,1]$ for $N=201$. The collocation points are the Legendre points and an oversampling factor of 2 was chosen. The last panel (d) is analogous to (c) but is based on Chebyshev extension with Chebyshev roots.}
\end{figure}

The setting can be generalized to other bases that satisfy a discrete orthogonality condition similar to \eqref{eq:discrete_orthogonality}, for example orthogonal polynomials with their roots as sampling points. We consider the approximation of a smooth function on $\left[-\tfrac12,\tfrac12\right]$, using Legendre polynomials up to degree $N-1$ which are orthogonal on $[-1,1]$. We sample the function in $M$ Legendre nodes in the subinterval $\left[-\tfrac12,\tfrac12\right]$. In order to ensure oversampling, these nodes are the roots of a Legendre polynomial of higher degree $L$, where $L$ is chosen such that the restricted point set has size $M = L/2 > N$. The resulting singular value patterns of $A$, $Z$ and $A-AZ^*A$ are shown in the first three panels of Figure \ref{fig:legendre_az}. The singular values of $A$ and $Z$ do not cluster around $1$ as in the case of Fourier extension. Yet, we can still isolate the plunge region that is present in $A$ by the matrix $A-AZ^*A$.

The discrete orthogonality condition of Legendre polynomials on the full grid of length $L$ is
\begin{equation}\label{eq:discrete_orthogonality_legendre}
 \sum_{l=1}^L w_l P_i(x_l) P_j(x_l) = h_i^2 \delta_{i-j}, \qquad 0 \leq i,j \leq L-1,
\end{equation}
where $w_l$ are the Gauss--Legendre weights associated with the roots $x_l$, and $h_i = \Vert P_i \Vert_{[-1,1]}$ is the norm of $P_i$. As in the Fourier case, there is a large $L \times L$ matrix $F$ with entries $F_{i,j} = P_{i-1}(x_{j-1})$. This matrix has inverse $F^{-1} = DF^*W$, where $W$ and $D$ are diagonal matrices with entries $w_i$ and $h_i^{-2}$ respectively. The discrete least squares matrix $A$ is a submatrix of $F$, with columns selected corresponding to the degrees of freedom ($0,\ldots,N-1$) and rows selected corresponding to the points in the subinterval $[-\frac12,\frac12]$. We choose $Z$ as the corresponding subblock of $(F^{-1})^* = W F D$. The size of the plunge region in this setting is not known in literature and is the topic of a separate study.

The same methodology can be applied using Chebyshev polynomials and Chebyshev roots. Here, due to the relation between trigonometric polynomials and Chebyshev polynomials, there is a fast matrix-vector product for $A$ and for $Z$ based on the discrete cosine transform. The singular value patterns of $A$ and $Z^*$ are not shown for this case, but they are similar to the Fourier extension case in Figure \ref{fig:fourier_az}. The singular value pattern of $A-AZ^*A$ is in the bottom right panel of Figure \ref{fig:legendre_az}. As in the Fourier extension case the rank of this matrix is $\mathcal O(\log(N))$.

Finally, the Chebyshev roots can also be replaced by the Chebyshev extremae, in which case a discrete orthogonality property based on Clenshaw--Curtis quadrature (as opposed to Gauss--Chebyshev) yields a suitable matrix $Z$. The different choice of discretization points corresponds to a different norm for the continuous approximation problem: we find the best approximation either with respect to the Legendre-weighted norm ($w(x)=1$, Clenshaw--Curtis) or the Chebyshev-weighted norm with  ($w(x) = \frac{1}{\sqrt{1-x^2}}$, Gauss--Chebyshev).

% !TEX root = AZalgorithmV3.tex
\section{Weighted linear combinations of bases}
\label{s:sumframes}

Let $\Phi = \{ \phi_k \}_{k=0}^\infty$ be a basis on $[-1,1]$. Here, we consider \emph{weighted sum frames} of the form
\begin{equation}\label{eq:sumframe}
 \Psi = w_1 \Phi \cup w_2 \Phi = \{ w_1(x) \phi_k(x) \}_{k=0}^\infty \cup \{ w_2(x) \phi_k(x) \}_{k=0}^\infty,
\end{equation}
where $w_1$ and $w_2$ are two functions of bounded variation on $[-1,1]$ such that
\begin{equation} \label{eq:sumframe_condition}
w_1(x)^2 + w_2(x)^2 > 0, \qquad \forall x \in [-1,1].
\end{equation}

\subsection{The canonical dual frame and $(A,Z)$}

Assume that the basis $\Phi$, truncated after $N$ terms to $\Phi_N$, has an associated least squares grid $\{ x_m \}_{m=1}^M$ of length $M > N$. The corresponding discretization matrix $A_\Phi \in \mathbb{C}^{M \times N}$ has entries $A_{m,n} = \phi_n(x_m)$. Assume also that the least squares system $A_\Phi x = b$ is solved by $Z_\Phi^* b$. In other words, we assume that the basis $\Phi$ has an $(A_\Phi,Z_\Phi)$ combination such that $Z_\Phi^*$ yields exact reconstruction for all functions in the span of $\Phi_N$.

If the basis $\Phi_N$ is a Fourier basis or Chebyshev polynomial basis, then a least squares problem for this basis is solved by a truncated FFT or discrete cosine transform --- in this case $A_\Phi$ and $Z_\Phi$ have efficient implementations.

The $A$ matrix for the weighted sum frame is given by
\begin{equation}\label{eq:sumframe_A}
 A = \begin{bmatrix}
 W_1 A_\Phi & W_2 A_\Phi
 \end{bmatrix}
\end{equation}
where $W_{1,2} \in \mathbb{C}^{M \times M}$ are diagonal matrices with entries $w_{\{1,2\}}(x_m)$ on the diagonal.

The $Z$ matrix follows from the construction of a dual frame. The canonical dual frame for \eqref{eq:sumframe} is given by (see also \cite[Example 3]{adcock2019frames}):
\begin{equation*} \label{eq:sumframe_dual}
 \tilde{\Phi} = \left\{ \frac{1}{|w_1(x)|^2+|w_2(x)|^2} \, \phi_k(x) \right\}_{k=0}^\infty.
\end{equation*}
Note that condition \eqref{eq:sumframe_condition} was imposed in order to guarantee the existence of this dual frame. We illustrate by example in what follows, that a suitable $Z$ matrix for the weighted sum frame is
\begin{equation}\label{eq:sumframe_Z}
 Z = \begin{bmatrix}
 W^\dagger{W_1} Z_\Phi & W^\dagger{W_2} Z_\Phi
 \end{bmatrix}
\end{equation}
where $W \in \mathbb{R}^{M \times M}$ is a diagonal matrix with entries
\[
W_{m,m} = |w_1(x_m)|^2+|w_2(x_m)|^2.
\]

\subsection{Weighted combination of bases}

Consider the function
\[
f(x,y) = \cos(2\pi(x+y)) + \sqrt{x^2+y^2} \sin(1+2\pi(x+y))
\]
on the square domain $[-1,1]^2$. We approximate $f$ with a weighed sum frame
\[
\Psi = \Phi \cup \sqrt{x^2+y^2} \Phi,
\]
where $\Phi$ is a tensor product of Chebyshev polynomials.

\begin{figure}
	\begin{center}
		\subfigure[Singular values of $A$]{
			\includegraphics[width=.4\textwidth]{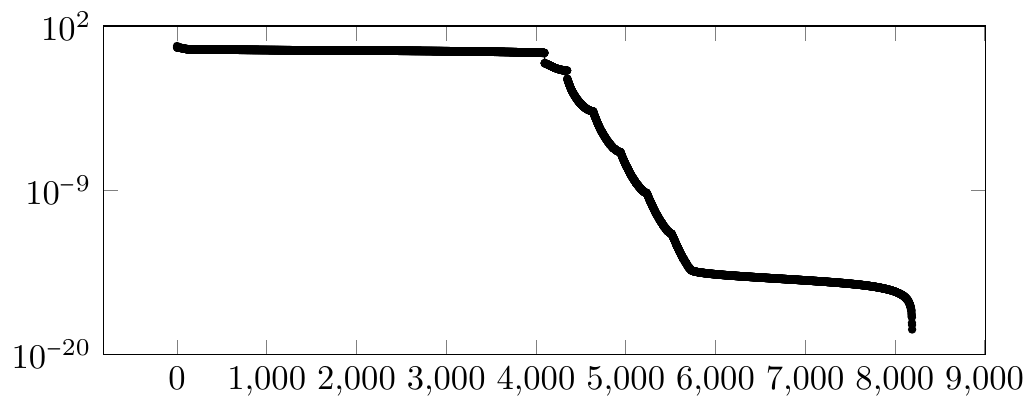}
		}
		\hspace{0.5cm}
		\subfigure[Singular values of $A-AZ^*A$]{
			\includegraphics[width=.4\textwidth]{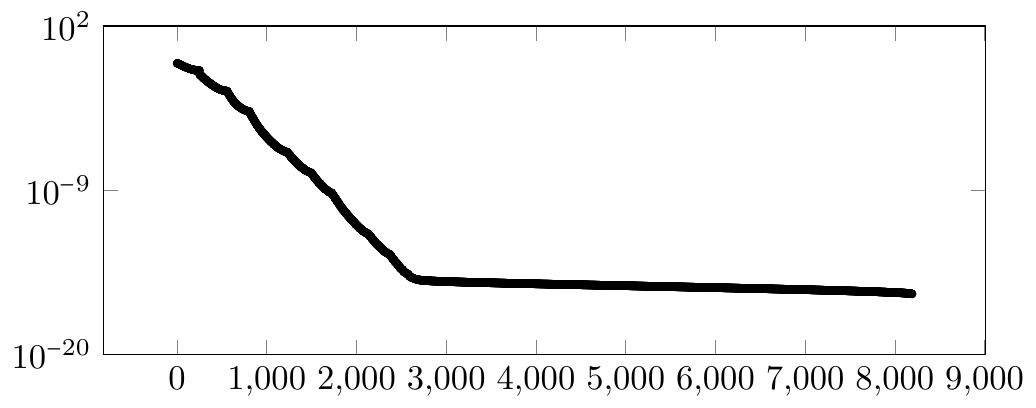}
		}\\
		\subfigure[Timings (in seconds)]{
			\includegraphics[width=.4\textwidth]{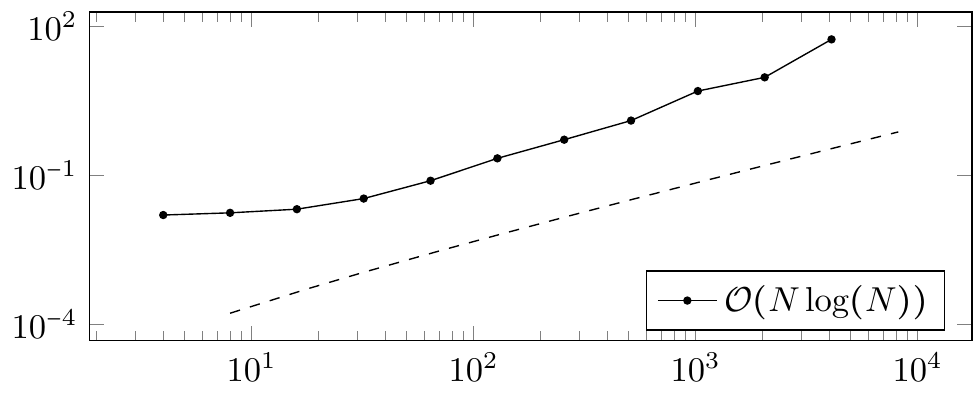}
		}
		\hspace{0.5cm}
		\subfigure[Approximant and error]{
			\includegraphics[width=.4\textwidth]{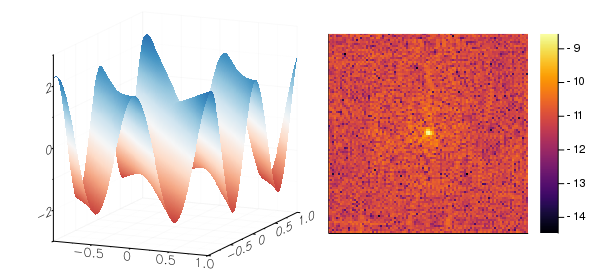}
		}
	\end{center}
	\caption{\label{fig:weightedsquare} Illustration of the AZ algorithm for approximation with $\Phi \cup \sqrt{x^2+y^2} \Phi$, where $\Phi$ is a $64\times64$ tensor product of Chebyshev polynomials on $[-1,1]^2$. The collocation points are the 2D cartesian product of a $128$-point Chebyshev grid. In the last panel: the approximation of $f(x,y)=\cos(2\pi(x+y)) + \sqrt{x^2+y^2}\sin(1+2\pi(x+y))$ and the $\log_{10}$ error of this approximation. }
\end{figure}

The function $f$ has a square root singularity at the origin. The approximation space is chosen such that it effectively captures the singularity. It is obvious that $f$ can be well approximated in the span of the weighted sum frame. However, standard algorithms cannot be used to compute such an approximation if we are only allowed to sample $f$ (rather than, for example, its smooth and singular part separately). We show that the AZ algorithm with the above choices of $A$ and $Z$ is an effective means to this. The results are shown in Figure \ref{fig:weightedsquare}.

The study of the size of the plunge region in this example, which demonstrates the computational complexity of the AZ algorithm, is the topic of a separate study. For example, see the forthcoming paper \cite{webb2019plunge}.

\subsection{Weighted combination of frames}

It is not necessary to assume that $\Phi$ is a basis. We observe that the formulas remain valid if instead $\Phi$ is a frame, with a suitable known $(A_\Phi,Z_\Phi)$ combination. In that case, the matrix $A$ given by \eqref{eq:sumframe_A} and $Z$ given by \eqref{eq:sumframe_Z} form a suitable $(A,Z)$ combination for the weighted sum frame.

\begin{figure}
	\begin{center}
		\subfigure[Singular values of $A$]{
			\includegraphics[width=.4\textwidth]{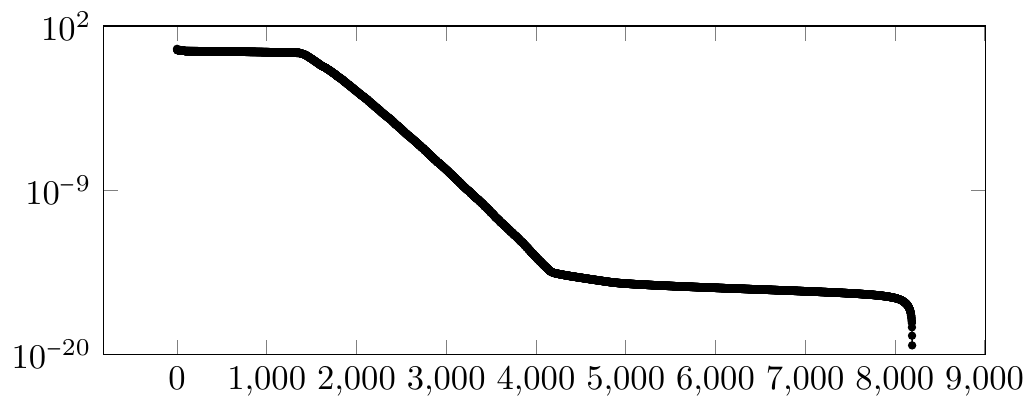}
		}
		\hspace{0.5cm}
		\subfigure[Singular values of $A-AZ^*A$]{
			\includegraphics[width=.4\textwidth]{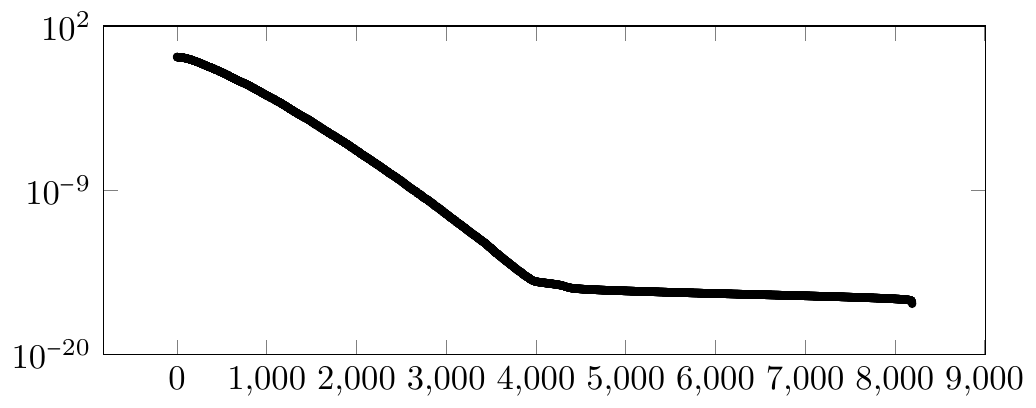}
		}\\
		\subfigure[Timings (in seconds)]{
			\includegraphics[width=.4\textwidth]{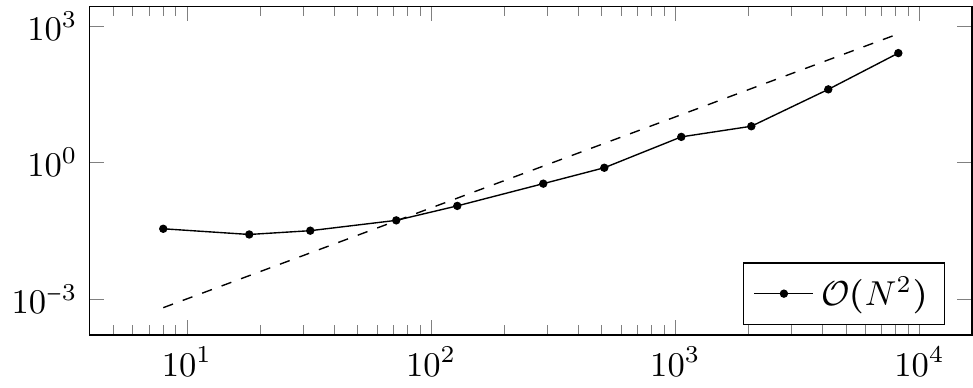}
		}
		\hspace{0.5cm}
		\subfigure[Approximant and error]{
			\includegraphics[width=.4\textwidth]{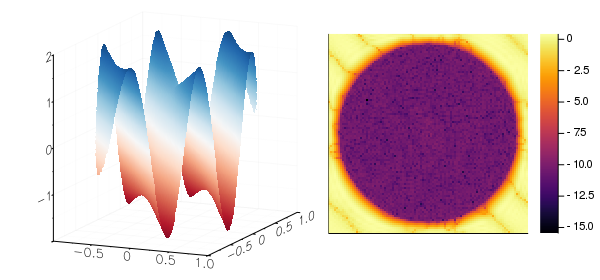}
		}
	\end{center}
	\caption{\label{fig:weightedcircle} Approximation with $\Phi \cup \sqrt{x^2+y^2} \Phi$ on a circular domain with center $(0,0)$ and radius $0.9$, where $\Phi$ is a $64\times64$ tensor product of Chebyshev polynomials on $[-1,1]^2$. The collocation points are the 2D cartesian product of a 128-point Chebyshev grid. In the last panel: the approximation of $f(x,y)=\cos(2\pi(x+y)) + \sqrt{x^2+y^2}\sin(1+2\pi(x+y))$ and the logarithmic error of this approximation. }
\end{figure}

We approximate the same singular function $f$ as in the previous section, but instead of using a rectangle we restrict samples of $f$ to a disk with center $(0,0)$ and radius $0.9$. Thus, we make a weighted linear combination of Chebyshev extension frames. A suitable $(A,Z)$ combination for Chebyshev extension was described in \S\ref{s:fourext}, and the experiment shows it can be used as part of a more complicated $(A,Z)$ combination for the weighted combination. This is shown in Figure \ref{fig:weightedcircle}.

% !TEX root = AZalgorithmV3.tex
\section{B-spline extension frames}
\label{s:splines}

In \S\ref{s:fourext} we considered extension frames using an orthogonal basis on a bounding box. In that section we emphasized the connection to continuous and discrete dual bases for the selection of $Z$. This connection is explored in more detail for extensions based on B-spline bases in \cite{coppe2019splines}. Since B-splines are compactly supported, the collocation matrix $A$ is highly sparse. B-splines are not orthogonal, but several different dual bases can be identified and it is shown in \cite{coppe2019splines} that each of these leads to a suitable $Z$ matrix. Some choices of dual bases lead to a sparse $Z$ and even a sparse matrix $A - AZ^*A$, with corresponding advantages for speed and efficiency. We refer to \cite{coppe2019splines} for the analysis and examples.

% !TEX root = AZalgorithmV3.tex
\section{Weighted least squares approximation}
\label{s:weighted}

In this final example we show how an existing efficient solver for a least squares approximation can be used to solve weighted variants of the problem efficiently as well. Thus, assume that a least squares problem is given in the form
\begin{equation}\label{eq:ls_A}
Ax=b,
\end{equation}
where $A \in \mathbb{C}^{M \times N}$ with $M > N$, and that a $Z$ matrix is known that approximates $A (A^*A)^{-1}$, i.e., such that $Z^*b$ solves the unweighted least squares problem. In this section we consider a Fourier series of length $N$, and a discrete least squares approximation to a given function based on $M > N$ equispaced samples. Owing to the continuous and discrete orthogonality properties of Fourier series, this problem is solved efficiently simply by computing the inverse FFT of length $M$ and truncating the result to a vector of length $N$.

Next, we consider a diagonal weight matrix $W \in \mathbb{C}^{M \times M}$ with positive entries $W_{i,i} = d_i > 0$, that associates a weight to each condition in the rows of $A$. The weighted least squares problem is
\begin{equation}\label{eq:ls_WA}
WAx=Wb.
\end{equation}
When the linear system is overdetermined, the weighted least squares solution differs from the unweighted one. The solution is given by
\[
 x = (A^* W^* W A)^{-1} A^* W^*W b.
\]
The solution of \eqref{eq:ls_A} minimizes $\Vert Ax-b\Vert$, the solution of \eqref{eq:ls_WA} minimizes $\Vert W(Ax-b)\Vert$.

With an efficient solver for \eqref{eq:ls_A} at hand, the simplest solution method to solve the weighted problem is to ignore the weight matrix and to return the solution to \eqref{eq:ls_A}. This may deviate from the weighted solution as follows,
\[
\min_{1 \leq i \leq M} d_i \, \Vert Ax-b \Vert \leq \Vert W(Ax-b)\Vert \leq \max_{1 \leq i \leq M} d_i \, \Vert Ax-b \Vert.
\]
The weighted and unweighted least squares problems may have radically different solutions if the ratio of the weights is large, and this is the setting we focus on. Assume w.l.o.g.~that the large ratio of weights is due to having a number of very small weights. We choose a threshold $\epsilon$ and we define $W_\epsilon$ as the diagonal matrix with entry $i$ equal to zero if $|d_i| < \epsilon$ and equal to $d_i$ otherwise. We define the matrix $\tilde{Z}$ in terms of the pseudoinverse of $W_\epsilon$,
\begin{equation*}\label{eq:ls_Ztilde}
 \tilde Z = W_\epsilon^\dagger Z.
\end{equation*}

We set out to solve $\tilde{A} x = \tilde{b}$, with $\tilde{A} = WA$ and $\tilde{b} = Wb$. With these choices, the AZ algorithm is given explicitly in Algorithm \ref{alg:AZweighted}.

\begin{algorithm}[h!]
	\caption{The AZ algorithm for weighted approximation}\label{alg:AZweighted}
	{\bf Input:} $A, Z \in \C^{M\times N}$, $W=\diag([d_1,\dots,d_M])$, with $d_i>0$,  $b\in\C^M$, $\epsilon > 0$ \\
	{\bf Output:} $x\in\C^N$ such that $WAx \approx Wb$
	\begin{algorithmic}[1]
		\State Solve $(I-WA(W_\epsilon^\dagger Z)^*)WAx_1 = (I-WA(W_\epsilon^\dagger Z)^*)Wb$
		\State $x_2 \gets Z^*W_\epsilon^\dagger W(b-Ax_1)$
		\State $x \gets x_1 + x_2$
	\end{algorithmic}
\end{algorithm}

One interpretation of the algorithm is as follows. The known $Z$ matrix is efficient, but it solves the wrong problem, namely the unweighted one. It is used in step 2. Step 1 is slow, but it solves the right problem, namely the weighted one. If $\epsilon=0$, then $W^\dagger_\epsilon = W^{-1}$. In this case, the system in the first step has rank zero and the problem is solved in step 2, yielding the solution to the unweighted problem: efficient, but possibly inaccurate. If on the other hand $\epsilon \gg 1$,  then $W^\dagger_\epsilon = 0$ and the problem is solved in step 1: accurate, but possibly slow. By varying $\epsilon$, one obtains a solution somewhere in between these two extreme cases. The tradeoff is between accuracy of the weighted problem and efficiency of the algorithm.

\begin{figure}
	\centering
	\includegraphics[width=.4\textwidth]{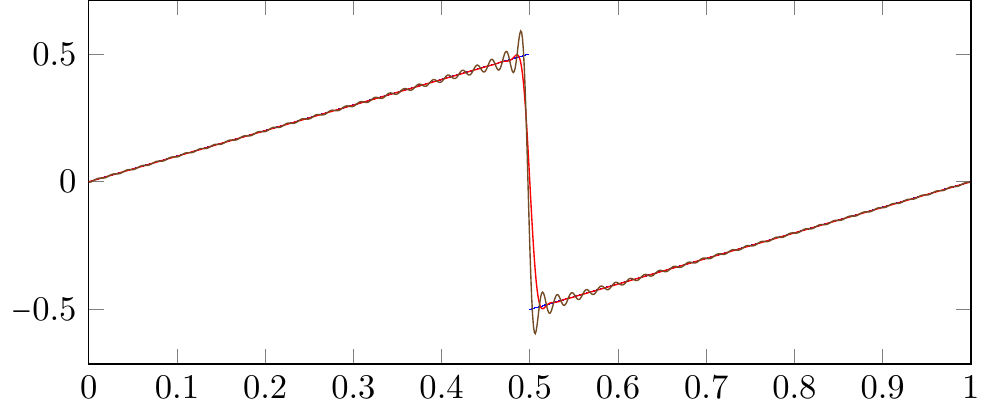}%
	\includegraphics[width=.4\textwidth]{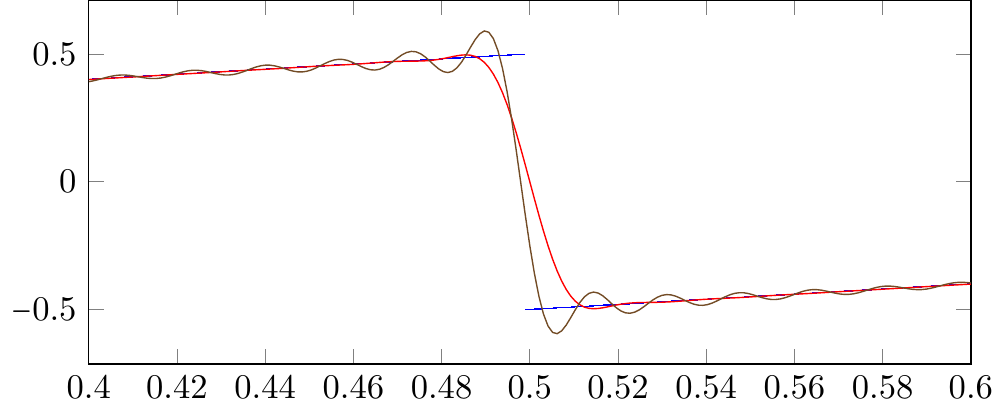}%
	\caption{\label{fig:weightedazfunction} Left panel: A function with a jump at $x=0.5$ (blue), the unweighted approximation (brown) containing the Gibbs phenomenon, and the weighted approximation (red). Right panel: Zoom of the left panel.
	}
\end{figure} 

For an example we return to function approximation: we use a Fourier series to approximate a piecewise smooth, yet discontinuous function. This is known to suffer the Gibbs phenomenon, resulting in an oscillatory overshoot at the discontinuities. In order to obtain a smooth approximation to the smooth parts of $f$, we may want to assign a very small weight to the function accuracy near the discontinuity. This is a simple alternative to other smoothing methods, such as spectral filtering techniques.

\begin{figure}
	\centering
	\includegraphics[width=.4\textwidth]{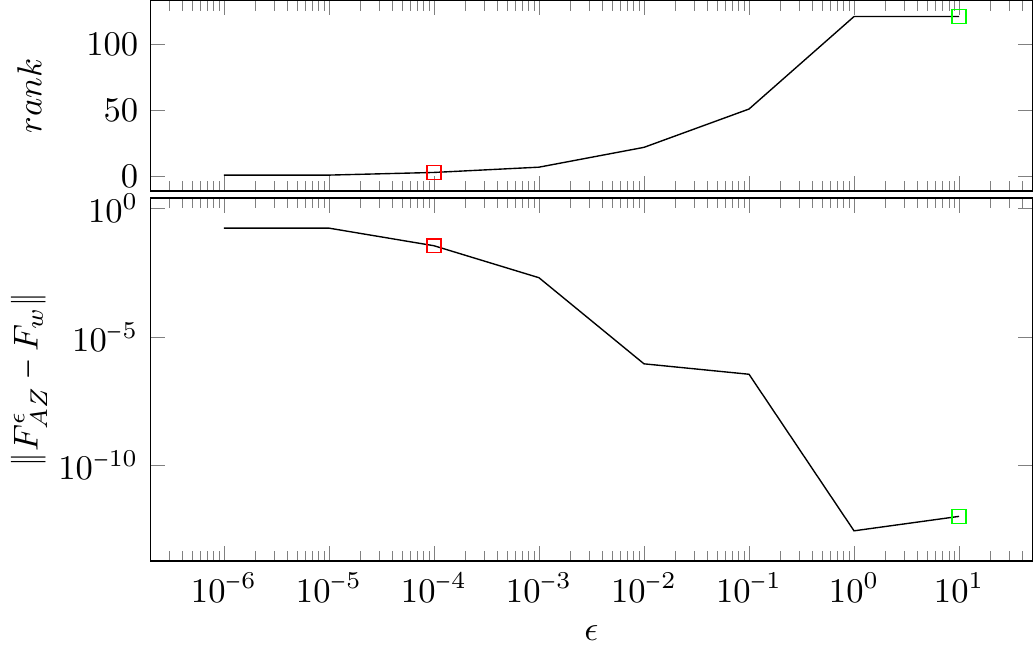}
	\includegraphics[width=.4\textwidth]{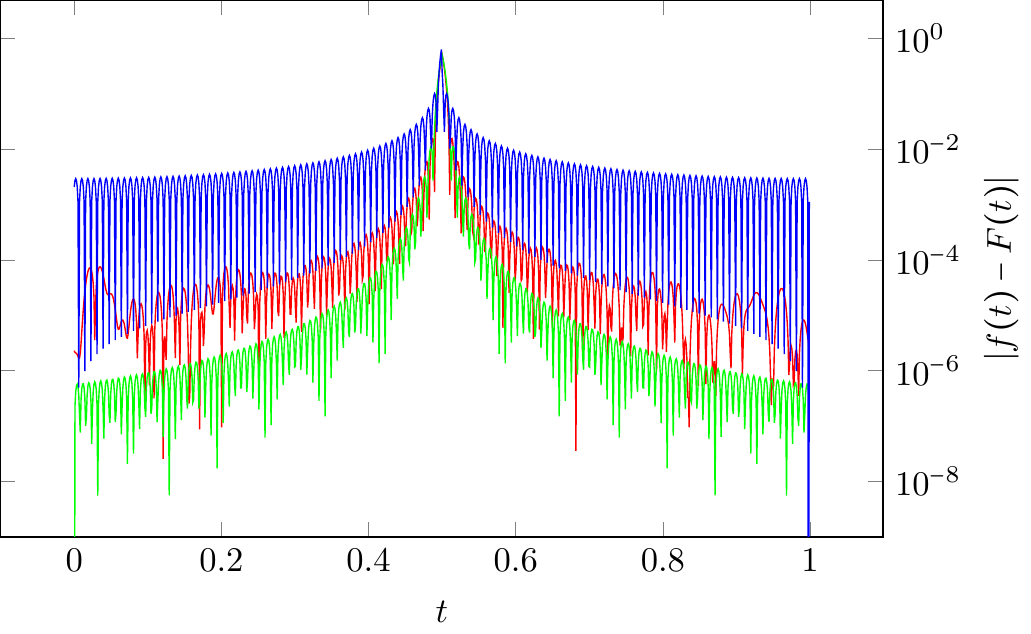}
	\caption{\label{fig:weightedaz} Left panel. Top: the rank of step 1 of the AZ algorithm. Bottom: the norm of the difference between the solution given by the AZ algorithm and the solution of $WAx=Wb$. Right panel. The difference between the weighed solution and the unweighted solution (blue), the error of the AZ solution with $\epsilon=10^{-4}$ (red), and of the AZ solution with $\epsilon=10$ (green).
	}
\end{figure}

The function shown in Fig.~\ref{fig:weightedazfunction} is periodic on $[0,1]$, with a jump at $0.5$. We use a Fourier series with $N=121$ terms on $[0,1]$. The discrete weights are based on sampling the function $w(x)=(x-0.5)^2$, which assigns greater weight to points further away from the discontinuity. We solve the weighted approximation problem with the AZ algorithm varying $\epsilon$ from $1e^{-6}$ to $1e^{+1}$. The results are shown in Fig.~\ref{fig:weightedaz}. For small $\epsilon$, the rank of the system in the first step is small, but the solution deviates substantially from the true weighted solution. For larger $\epsilon$, the rank of the system in step 1 increases, but the AZ algorithm approximates the true solution better. Fig.~\ref{fig:weightedazfunction} shows that the AZ algorithm produces a non-oscillatory, smooth approximation to $f$ while Fig.~\ref{fig:weightedaz} shows that this approximation is highly accurate away from the jump.

\section{Concluding remarks}
\label{s:conclusion}

We introduced a three step algorithm and showed that it can be used to efficiently solve various problems of interest in function approximation. The examples were chosen such that the matrix $Z$ could be devised analytically. However, as the examples show, it is often possible to generate $(A,Z)$ combinations for complicated approximation problems from $(A,Z)$ combinations of simpler subproblems.

The efficiency of the AZ algorithm hinges on the numerical rank of the system that is solved in the first step. Plots of the rapidly decaying singular values are provided for each example. Theoretical bounds on these numerical ranks is a topic of ongoing research, and beyond the scope of this paper. See \cite{matthysen2016fastfe,matthysen2017fastfe2d} for discussion of rank for the univariate and multivariate Fourier extension problem and see \cite{webb2019plunge} for discussion of rank for weighted combinations of certain univariate bases.

In our examples above, we have used the AZ algorithm in combination with a randomized low-rank SVD solver for step 1. Iterative solvers such as LSQR and LSMR \cite{Paige1982,Fong2010} can also be used for the first step of the AZ algorithm. Experiments indicate that they do not typically yield high accuracy, due to the ill-conditioning of the system. However, if only a few digits of accuracy are required, they can be several times more efficient than using a direct solver in step 1.

\bibliographystyle{siamplain}
\bibliography{AZbib}

\begin{thebibliography}{10}

\bibitem{adcock2018fna2}
{\sc B.~Adcock and D.~Huybrechs}, {\em Frames and numerical approximation {II}:
  generalized sampling}, Tech. Report arXiv:1802.01950, 2018.

\bibitem{adcock2019frames}
{\sc B.~Adcock and D.~Huybrechs}, {\em Frames and numerical approximation},
  SIAM Rev., 61 (2019), pp.~443--473.

\bibitem{boyd2002fourierextension}
{\sc J.~P. Boyd}, {\em A comparison of numerical algorithms for {F}ourier
  extension of the first, second, and third kinds}, J. Comput. Phys., 178
  (2002), pp.~118--160.

\bibitem{bruno2007continuation}
{\sc O.~P. Bruno, Y.~Han, and M.~M. Pohlman}, {\em Accurate, high-order
  representation of complex three-dimensional surfaces via {F}ourier
  continuation analysis}, J. Comput. Phys., 227 (2007), pp.~1094--1125.

\bibitem{christensen2008framesandbases}
{\sc O.~Christensen}, {\em Frames and Bases: an Introductory course}, Springer,
  Basel, 2008.

\bibitem{coppe2019splines}
{\sc V.~Copp\'e and D.~Huybrechs}, {\em Efficient function approximation on
  general bounded domains using splines on a cartesian grid}, Tech. Report
  arXiv:1911.07894, KU Leuven, 2019.

\bibitem{edelman1999futurefft}
{\sc A.~Edelman, P.~McCorquodale, and S.~Toledo}, {\em The future {F}ast
  {F}ourier {T}ransform?}, SIAM J. Sci. Comput., 20 (1999), pp.~1094--1114.

\bibitem{Fong2010}
{\sc D.~Fong and M.~Saunders}, {\em {LSMR}: An iterative algorithm for sparse
  least-squares problems}, SIAM J. Sci. Comput., 33 (2010), pp.~2950--2971.

\bibitem{Golub2013}
{\sc G.~H. Golub and C.~F. van Loan}, {\em Matrix Computations}, JHU Press,
  Baltimore, {F}ourth~ed., 2013.

\bibitem{halko2011finding}
{\sc N.~Halko, P.-G. Martinsson, and J.~A. Tropp}, {\em Finding structure with
  randomness: Probabilistic algorithms for constructing approximate matrix
  decompositions}, SIAM Rev., 53 (2011), pp.~217--288.

\bibitem{huybrechs2010fourierextension}
{\sc D.~Huybrechs}, {\em On the {F}ourier extension of non-periodic functions},
  SIAM J. Numer. Anal., 47 (2010), pp.~4326--4355.

\bibitem{lawson1996leastsquares}
{\sc C.~L. Lawson and R.~J. Hanson}, {\em Solving least squares problems},
  Classics in Applied Mathematics, SIAM, Philadelphia, 1996.

\bibitem{lyon2011fastcontinuation}
{\sc M.~Lyon}, {\em A fast algorithm for {F}ourier continuation}, SIAM J. Sci.
  Comput., 33 (2011), pp.~3241--3260.

\bibitem{matthysen2016fastfe}
{\sc R.~Matthysen and D.~Huybrechs}, {\em Fast algorithms for the computation
  of {F}ourier extensions of arbitrary length}, SIAM J. Sci. Comput., 38
  (2016), pp.~A899--A922.

\bibitem{matthysen2017fastfe2d}
{\sc R.~Matthysen and D.~Huybrechs}, {\em Function approximation on arbitrary
  domains using fourier frames}, SIAM J. Numer. Anal., 56 (2018),
  pp.~1360--1385.

\bibitem{Paige1982}
{\sc C.~C. Paige and M.~A. Saunders}, {\em {LSQR}: An algorithm for sparse
  linear equations and sparse least squares}, ACM Trans. on Math. Software, 8
  (1982), pp.~43--71.

\bibitem{slepian1978prolateV}
{\sc D.~Slepian}, {\em Prolate spheroidal wave functions, {F}ourier analysis,
  and uncertainty {V}: the discrete case}, Bell Labs Tech. J., 57 (1978),
  pp.~1371--1430.

\bibitem{webb2019plunge}
{\sc M.~Webb}, {\em The plunge region in approximation by singularly enriched
  trigonometric frames}, Tech. Report to appear, 2019.

\end{thebibliography}

\end{document}